\documentclass[12pt]{article}
\usepackage{amsmath,amsfonts,amssymb,xspace}
\usepackage{verbatim}
\usepackage{mathtools}
\usepackage{enumitem} 
\usepackage{color}
\definecolor{red}{rgb}{.8,0,0}

\usepackage{graphicx}
\usepackage{wrapfig}
\usepackage{tikz}
\usepackage{amsmath}
\usepackage{url}
\usepackage[mathlines]{lineno}
\usepackage{color}

\newtheorem{theorem}{Theorem}[section]
\newtheorem{lemma}[theorem]{Lemma}
\newtheorem{prop}[theorem]{Proposition}
\newtheorem{cor}[theorem]{Corollary}
\newtheorem{remark}[theorem]{Remark}
\newtheorem{definition}[theorem]{Definition}
\newtheorem{example}[theorem]{Example}

\newenvironment{proof}{\noindent{Proof:}}{\hfill\qed\smallskip}
\newcommand{\qed}{\quad\rule{1.5ex}{1.5ex}}

\newcommand{\mbb}{\mathbb}

\newcommand{\mc}{\mathcal}

\title{Partial geometric designs having circulant concurrence matrices}
\author{Sung-Yell Song\footnote{sysong@iastate.edu, Department of Mathematics, Iowa State University, Ames, IA 50011, U. S. A.} \ and Theodore Tranel\footnote{teddyt@iastate.edu, Department of Mathematics, Iowa State University, Ames, IA 50011, U. S. A.}}
\date{12-29-2021}

\begin{document} %
\pagenumbering{arabic} \setcounter{page}{1}

\maketitle

\begin{abstract}  
We survey partial geometric designs and investigate their concurrences of points. The concurrence matrix of a design, which encodes the concurrences of pairs of points, can be used in the classification of designs in some extent. An ordinary 2-$(v,k,\lambda)$ design has concurrence $\lambda$ for any pair of distinct points, and its concurrence matrix is circulant. A partial geometry has two concurrences $1$ and $0$ and a transversal design TD$_{\lambda}(k, u)$  has two concurrences $\lambda$ and $0$. It is also known that the concurrence matrix of a partial geometric design can have at most three distinct eigenvalues, all of which are non-negative integers. In this paper, we show the existence of other partial geometric designs having two or three distinct concurrences, and investigate which symmetric circulant matrices are realized as the concurrence matrices of partial geometric designs. We collect known sources of partial geometric designs and study their structural characteristics and construction methods. We then give a list of feasible parameter sets for partial geometric designs of order up to 12 each of which has a circulant concurrence matrix. We also consider the combinatorial properties and constructions of partial geometric designs satisfying these parameter sets.\footnote{This work is part of the Ph.D. dissertation of the second author \cite{Tr}.}
\end{abstract}
\bigskip

{\flushleft Keywords: association scheme, finite incidence structure, $t$-design, partial geometry, partial geometric difference set, special partially balanced incomplete block design. strongly regular graph.
\bigskip

}

\newpage
\tableofcontents 
\newpage

\section{Introduction}
Partial geometric designs were introduced by Bose, Shrikhande and Singhi \cite{BSS} in their study of an application to the embedding theorem for balanced incomplete block designs (Hall and Connor \cite{CH}).  Upon introduction, partial geometric designs could be regarded as a generalization of a partial geometry, which was introduced by Bose \cite{Bose}.  In \cite{N}, Neumaier generalized the concept of a partial geometric design in his study of $t\frac{1}{2}$-designs and called a partial geometric design a  $1\frac{1}{2}$-design. He laid the foundation in the study of this subject providing important basic properties of the partial geometric designs and their links to 2-designs, partial geometries, transversal designs and complete bipartite graphs. There have been extensive studies on properties and applications of partial geometric designs since then. In particular, Neumaier showed that the concurrence matrices of the partial geometric designs, that are not $2$-designs, have exactly three distinct eigenvalues with one of the eigenvalues being $0$.  We observe that the concurrence matrices of many partial geometric designs including 2-designs and transversal designs are circulant. Motivated by this observation, we intend to collect as many partial geometric designs whose concurrence matrices are circulant as possible. We begin our investigation by finding all the symmetric circulant matrices that have three prescribed integral eigenvalues.  Then we determine whether each of such circulant matrix can be realized as the concurrence matrix of a PGD. We give a list of our findings and compare with the classification results of van Dam and Spence reported in \cite{vDS}. We use the results from \cite{N} as a large part of the foundation for our investigation on the classification of PGDs based on the structure of their concurrence matrices.

The organization of the paper is as the following. In Section 2, we briefly recall a few finite incidence structures, and then review the basic properties of partial geometric designs. In Section 3, we survey the sources of partial geometric designs and investigate the connections between the sources. The sources that will be discussed include transversal designs, partial geometric difference sets, strongly regular graphs, association schemes, partial geometries, (special) partially balanced incomplete block designs, affine resolvable designs, and classical finite geometries. In Section 4, we discuss some characterization results of partial geometric designs in terms of concurrences of points, including the work of Qu, Lei and Shan \cite{LQS}. In Section 5, we give a complete list of parameter sets for partial geometric designs of order up to 12 whose concurrence matrices can be circulant, and then construct partial geometric designs for each of the feasible parameter sets.

\section{Preliminaries}

We briefly recall necessary terms, set the notation, and review the basic facts that we will use throughout the paper.

\subsection{Finite incidence structures}

A finite incidence structure $S$ consists of a finite set $P$ of points,  a finite set $\mathcal{B}$ of blocks (or lines), and an incidence relation between the points and blocks. We denote $S$ $=$ $\left(P, \mathcal{B}, {\bf \mathcal{I}}\right)$ together with incidence relation $\mathcal{I} \subseteq P \times \mathcal{B}$ as the set of flags, i.e., the incident point-block pairs. If two different members, $B_1$ and $B_2$ of $\mathcal{B}$ are incident with the same set of points of $P$, we call these \emph{repeated blocks}.  A \emph{simple} finite incidence structure is one that does not contain repeated blocks. In a simple incidence structure, each block will be identified with the set of points that are incident with the given block. In this case, we can just consider blocks as subsets of $P$, replace the flag notation  $\left(p,B\right) \in \mathcal{I}$ by $p \in B$, and denote $S = \left(P, \mathcal{B}, \mathcal{I}\right)$ by $S = \left(P, \mathcal{B}\right)$ where the natural incidence relation is assumed. Throughout, we use $v$ and $b$ to denote the cardinality of $P$ and $\mathcal{B}$, respectively, unless otherwise specified.

For a finite incidence structure $S$ $=$ $\left(P, \mathcal{B}, {\bf \mathcal{I}}\right)$ where $P=\{p_1, p_2, \hdots, p_v\}$ and $\mathcal{B} = \{B_1, B_2, \hdots, B_b\}$,  the incidence matrix $N$ of $S$ is a $v \times b \,\, \{0, 1\}$-matrix:
$$N = [n_{ij}] \hspace{.75cm} \text{defined by the rule} \hspace{.5cm} n_{ij} =   \begin{cases} 
      1 & \text{if} \,\,\, \left(p_i, B_j\right) \in \mathcal{I} \\
      0 & \text{otherwise.} 
   \end{cases}
$$
The concurrence matrix of $S$ is then defined as $NN^T$ $($where $N^T$ is the transpose of $N)$; that is, the $(i, j)$-entry of $NN^T$ is $\lambda_{p_ip_j}:=|\{B\in\mathcal{B}: p_i, p_j\in B\}|$, which will be referred to the concurrence of the pair of points $p_i$ and $p_j$.
Throughout, let $I (= I_n)$ denote the ($n \times n$) identity matrix and let $J$ denote the all-ones matrix of an appropriate size. 

Given a finite incidence structure $S$ $=$ $\left(P, \mathcal{B}, {\bf \mathcal{I}}\right)$ with its incidence matrix $N$,  the \textit{dual} of $S$ is $\overline{S} = \left(\overline{P}, \overline{\mathcal{B}}, \overline{\mathcal{I}}\right)$ where $\overline{P} = \mathcal{B}$, $\overline{\mathcal{B}} = P$, and $\left(p,B\right) \in \mathcal{I}$ if and only if $\left(B,p\right) \in \overline{\mathcal{I}}$.  The dual of $S$ has incidence matrix $N^T$.
The \textit{complementary structure} of $S$ is $S^C = \left(P, \mathcal{B}^C, \mathcal{J}\right)$ where $\mathcal{B}^C = \{B^C = P - B : B \in \mathcal{B}\}$ and $\mathcal{J} = P \times \mathcal{B} - \mathcal{I}$, i.e. $\left(p,B^C\right) \in \mathcal{J}$ in $S^C$ if and only if $\left(p,B\right) \notin \mathcal{I}$ in $S$. The complement of $S$ has incidence matrix $J - N$.

Two incidence structures  $S=$ $\left(P, \mathcal{B}, {\bf \mathcal{I}}\right)$ and  $S'=$ $\left(P', \mathcal{B}', {\bf \mathcal{I}'}\right)$ are \textit{isomorphic} if there exists a bijection $f: P \rightarrow P'$ such that $\{f(B) : B \in \mathcal{B} \} = \{B': B' \in \mathcal{B}' \}$. In terms of their incidence matrices $N$ and $N'$, $S$ and $S'$ are isomorphic if and only if there exist permutation matrices $P$ and $Q$ such that $N'=PNQ$.

A finite incidence structure $S = \left(P, \mathcal{B}\right)$ with $|P| = v$, $|\mathcal{B}| = b$, is called
a \textit{tactical configuration} with parameters $\left(v, b, k, r\right)$ if  each block consists of exactly $k$ points, and each point belongs to exactly $r$ blocks.
A tactical configuration $S = (P, \mathcal{B})$ with $\left(v, b, k, r\right)$ satisfies 
$vr=bk$,\  $NJ = kJ$, and $JN = rJ$.

Given positive integers $v, k, t$, and $\lambda$ ($v > k > t$),  a \textit{$t$-design}, denoted by $t$-$\left(v,k,\lambda \right)$, is a tactical configuration satisfying  any $t$-subset $T$ of $P$ is contained in exactly $\lambda$ blocks. In this case, we have: $b = \lambda {v \choose t}/{k \choose t}$.
In this context,  a tactical configuration with parameters $\left(v, b, k, r \right)$ is a $1$-$\left(v, k, r \right)$ design. It is shown that a $t$-design is also a $t'$-design for $t' = 1, 2, \hdots, \left(t-1\right)$. With $t=2$, a $2$-$\left(v, k, \lambda \right)$ design is also commonly referred to as a \emph{balanced incomplete block design} (BIBD).
For the incidence matrix $N$ of a $2$-$\left(v, k, \lambda \right)$ design, we have $NN^T = rI + \lambda (J-I)$.

\subsection{Properties of partial geometric designs}
The partial geometric designs are introduced by Bose, Shrikhande, and Singhi \cite{BSS}. Most of what we discuss in this section is due to Neumaier \cite{N}. 

\begin{definition}
A partial geometric design
with parameters $\left(v, b, k, r; \alpha, \beta \right)$ is a tactical configuration $\left(P,\mathcal{B}\right)$ with parameters $\left(v, b, k, r\right)$ such that for any  $a \in P$ and any  $B \in \mathcal{B}$,
$$s\left(a,B\right) = \begin{cases} 
      \alpha & \text{if} \,\,\, a \notin B \\
      \beta & \text{if} \,\,\, a \in B 
   \end{cases}$$
where 
$$s\left(a,B\right): = |\{\left(x,A\right) \in P \times \mathcal{B}: x \in A, x \in B, a \in A\}|.$$
\end{definition}
 
We will refer to a partial geometric design, simply as a PGD or a PGD$(v, b, k, r;\alpha, \beta)$ from now on.
For reasons, by following Neumaier \cite{N}, which will later become clear, we let\footnote{In \cite{N}, the number $s\left(a, B\right)$ counts the number of flags $\left(x, A\right)$ satisfying the conditions $x\in B-\{a\}$, $A\ni a$ and  $A\neq B$; hence, $\beta$ is $k+r-1$ less than our $\beta$. As a consequence, $n=k+r-1+\beta-\alpha$ in \cite{N}.}
\begin{equation}\label{n}
    n = \beta-\alpha
\end{equation}

\begin{prop}{(\cite[Sec. 3.3]{N})}\label{parameters}
Let $\left(P,\mathcal{B}\right)$ 
be a PGD$(v, b, k, r; \alpha, \beta)$ and $n=\beta-\alpha$.
\begin{enumerate}
    \item[$1.$] \qquad
$\left(v-k\right)\alpha + k\beta = k^2r$
    \item[$2.$] \qquad
     $v=\frac{k\left(kr-n\right)}{\alpha},\ \ b=\frac{r\left(kr-n\right)}{\alpha},\ \  \beta = \alpha + n, \ \ 
     k+r \leq n + \alpha + 1 \leq kr.$
\end{enumerate}
\end{prop}
\begin{proof}
    1. It follows from counting the number of triples $\left(x, y, A\right)\in P\times P\times \mathcal{B}$ satisfying $x, y\in A$ and $y\in B$ for a fixed block $B$ in two ways. \\
    2. The relations follow from Part 1, and $vr=bk$, $ v>k>0, \,\, \alpha >0, \,\, \beta\ge k+r-1$.
\end{proof}

The following lemma will be used throughout. 

\begin{lemma}\label{conc} Let $\mathcal{D}=\left(P, \mathcal{B}\right)$ be a tactical configuration with parameters $\left(v, b, k, r\right)$. Let  $N$ be its incidence matrix, $\lambda_{xy} = \left |\{B\in\mathcal{B}: x, y\in B\}\right |$, and
 $$s\left(a,B\right): = |\{\left(x,A\right) \in P \times \mathcal{B}: x \in A, x \in B, a \in A\}|.$$
 Then we have the following.

\begin{enumerate}
\item[$(a)$] $\forall x, y\in P$, \[\lambda_{xy}  = \left [NN^T\right ]_{xy} \qquad (\mbox{and } \ \ \lambda_{xy} = r \mbox{ if } x=y).\]

\item[$(b)$] $\forall x\in P$, \[\sum\limits_{y\in P}\lambda_{xy}=r+\sum\limits_{y\in P-\{x\}}\lambda_{xy}= rk, \qquad (\mbox{equivalently, } NN^TJ = rkJ=JNN^T).\] 

\item[$(c)$] $\forall x\in P, \ B\in\mathcal{B},$ \[s\left(x, B\right) =\sum\limits_{y\in B}\lambda_{xy}= \left [NN^TN\right ]_{xB}.\]

\item[$(d)$] $\forall x, y\in P$, \[\sum\limits_{z\in P} \lambda_{xz}\lambda_{yz} =\left [ \left ( NN^T\right )^2\right ]_{xy} =\left [ \left (NN^TN\right)N^T\right ]_{xy} =\sum\limits_{B: B\ni y}s\left(x, B\right).\]
\end{enumerate}
 \end{lemma}
 
 \begin{proof} Straightforward from the definitions of $N$ and $s\left(x, B\right)$ and by routine combinatorial counting of flags of the incidence structure. For the details, we refer the reader to \cite{N}.
 \end{proof}

 \begin{lemma} \label{conc1} Suppose a tactical configuration $\mathcal{D}$ with $(v, b, k, r)$ is a PGD. Then there exist $\alpha, \beta \in \mathbb{N}\cup \{0\}$, such that: 
 \begin{enumerate}
 \item[] 
\[(i) \quad\qquad NN^TN=\beta N+\alpha \left(J-N\right)\quad \mbox{and} \quad
s\left(x, B\right)=\left \{\begin{array}{ll} \alpha &\mbox{if } x\notin B; \\ \beta & \mbox{if } x\in B.\end{array}\right. \]
 \item[] 
\[(ii)\qquad\qquad \sum\limits_{z\in P} \lambda_{xz}\lambda_{yz} =[\left(nN +\alpha J\right)N^T]_{xy}= \left \{\begin{array}{cl} \beta r & \mbox{if } x=y\\  
 n\lambda_{xy}+\alpha r& \mbox{if } x\neq y.\end{array}\right .\]
  \end{enumerate}
 \end{lemma}

\begin{proof} Proof is straightforward. For instance the second part follows from 
 \[  \sum\limits_{B: B\ni y}s\left(x, B\right) = \sum\limits_{B: B\ni y, x}s\left(x, B\right)+\sum\limits_{B: B\ni y,  x\notin B}s\left(x, B\right) = \beta \lambda_{xy} + \alpha \left(r-\lambda_{xy}\right).\]
 \end{proof}
 
 \begin{remark}
 We note that a 2-$\left(v, k, \lambda\right)$ design is a PGD with parameters \[\left(v,\ \  \frac{\lambda v\left(v-1\right)}{k-1},\ \ k,\ \ \frac{\lambda \left(v-1\right)}{k-1};\ \ \lambda k, \ \ \lambda \left(k-1\right)+r\right)\quad \mbox{and }\ \ n=r-\lambda = \frac{\lambda \left(v-k\right)}{k-1}.\] 
The concurrence matrix $NN^T=\left(r-\lambda \right)I+\lambda J$ of this design has two eigenvalues, $kr$ and $n=r-\lambda$ with multiplicities 1 and $v-1$, respectively. 
\end{remark}
Neumaier in \cite{N} showed that the concurrence matrix of a PGD, which is not a 2-design, has exactly three distinct eigenvalues including $0$ as in the following.

 \begin{lemma} \label{Spec} \cite[Sec. 3.12]{N} Let $N$ be the incidence matrix of a PGD$\left(v, b, k, r;\alpha, \beta\right)$. Then
 $NN^T$ has the eigenvalues $kr$, $n$, $0$ with multiplicities 1, $\sigma$, and $v-1-\sigma$, respectively, where $\sigma=r\left(v-k\right)/n$.
 In particular, $n>0$ and $r\left(v-k\right)\equiv 0\pmod{n}$.
 \end{lemma}
  
In algebraic graph theory, graphs with few distinct eigenvalues have been studied by numerous authors (cf. \cite{vDS} and the references in there). Van Dam and Spence noticed that there are regular bipartite graphs associated with nontrivial PGDs having four or five eigenvalues. This graph will be called the \emph{incidence graph} of the PGD.

\begin{definition} Let $N$ be the incidence matrix of a PGD with parameters $(v, b, k, r;\alpha, \beta)$. Then the regular bipartite graph whose adjacency matrix $A$ is given by 
\[A=\left [\begin{array}{cc} 0 & N\\ N^T & 0\end{array}\right ]\]
is called the incidence graph of the PGD.
\end{definition}

\begin{prop}\cite{vDS} For a given PGD$\left(v, b, k, r; \alpha, \beta \right)$ with its incidence matrix $N$, suppose the spectrum of $NN^T$ is $[kr^1,\ \ n^{\sigma}, \ \ 0^{v-1-\sigma}]$. Then the adjacency matrix $A$ of the incidence graph has eigenvalues $\pm \sqrt{kr}, \pm\sqrt{n}$, and possibly 0. 
\end{prop}

In regards to the possible concurrences for a PGD which is not a 2-design, we have the following observation from \cite[Sec. 3.14.]{N}. 
 \begin{lemma} Let $\left(P,\mathcal{B}\right)$ be a PGD which is not a 2-design and satisfies $\lambda_{xy}\in\{\lambda_1, \lambda_2\}$ for all $x\neq y$ and $\lambda_1>\lambda_2$. 
 Let $x\in P$ be a point. Suppose that there exist $m$ points $y\neq x$ with $\lambda_{xy}=\lambda_1$, and $v-1-m$ points $y\neq x$ with $\lambda_{xy}=\lambda_2$.  Then $0<m<v-1$ and 
\[m=\frac{r(k-1)-(v-1)\lambda_2}{\lambda_1-\lambda_2}.\]
 \end{lemma}

In what follows, we will mainly focus on PGDs with block size at least 3, and assume that $\alpha > 0, \,\,\, 3 \leq k \leq v-3, \,\,\, 3 \leq r \leq b-3$. We call these partial geometric designs \emph{proper} and others improper. 
We note that 
both the dual structure and the complementary structure of a PGD are PGDs as well.
\begin{lemma} Let $\mathcal{D}=\left(P, \mathcal{B}\right)$ be a proper PGD$\left(v, b, k, r; \alpha, \beta\right)$. Then we have the following parameters.
\begin{enumerate}
\item[$(1)$] The dual design of $\mathcal{D}$ is a PGD$\left(b, v, r, k;\alpha, \beta\right)$.
\item[$(2)$] The complementary design $\overline{\mathcal{D}}=\left(P, \overline{\mathcal{B}}\right)$ is a PGD$\left(v, b, v-k, b-r; \overline{\alpha}, \overline{\beta}\right)$ where
\[\overline{\alpha} = vb+3kr-bk-2vr-\beta,  \qquad \overline{\beta}= vb+3kr-bk-2vr-\alpha.\]
\end{enumerate}
\end{lemma}
\begin{proof} For (2), it is straightforward to verify that 
$\overline{N}\ \overline{N}^T\overline{N} = \overline{\beta}\ \overline{N} + \overline{\alpha}(J-\overline{N}).$
\end{proof}

\section{Known sources of partial geometric designs}

Partial geometric designs come from many combinatorial and geometric structures (cf. \cite{Bo, BSS, N}). In \cite{vDS}, by using an aid of the computer, van Dam and Spence classified and gave a complete list of small PGDs (having the sum of the numbers of points and blocks less than 36). In \cite{O}, Olmez introduced the notion of partial geometric difference sets to study the existence and construction of PGDs using group rings and group characters. There are more recent results on the constructions of PGDs from other incidence structures including some constructions of PGDs using totally isotropic subspaces of finite classical geometries over finite fields by Feng, Zhao and Zeng (cf. \cite{FZ, NO, NOS, NOS1, DO} and the references in it). 
We survey the known examples and constructions of PGDs some of which will be used later.

 
 \subsection{Transversal designs}
\begin{definition}
A \textit{transversal design}, denoted by TD$_\lambda(k, u)$, is a triple $\left(P, \mathcal{G}, \mathcal{B}\right)$ consisting of a point set $P$, a partition $\mathcal{G}$ of $P$ into $k$ classes each of size $u$, and a collection $\mathcal{B}$ of subsets of $P$ such that every block contains exactly one point from each class, and 
 every pair of points from different classes occur together in exactly $\lambda$ blocks.
 \end{definition}
It is clear that transversal design TD$_\lambda(k, u)$ is a tactical configuration with parameters $\left(v, b, k, r\right)  = \left(ku, \lambda u^2, k, \lambda u\right)$. The dual of a transversal design with $\lambda = 1$ is known as a \emph{net}.
Moreover, a TD$_{\lambda}(k,u)$  is a PGD with  parameters $\alpha=\lambda(k-1)$ and $\beta=\lambda(k-1) + \lambda u$.
The concurrence of two distinct points $x$ and $y$ is given by
\[\lambda_{xy} = \left \{ \begin{array}{ll} 0 & \mbox{if } x, y \mbox{ belong to the same group},\\ \lambda & \mbox{otherwise}.\end{array}\right .\]

A particular set of transversal designs comes from the vertex-edge incidence of the complete bipartite graph $K_{m,m}$. 
\begin{example} Given $m>2$, if we take the vertex set of the complete bipartite graph $K_{m,m}$ on $2m$ vertices as the point set and the set of edges as the block set, the graph becomes a PGD with parameters \[v=2m, \ \ b=m^2,\ \  k=2, \ \ r=m, \ \ \alpha =1, \ \ \beta = m+1.\] This PGD, referred to as a $K_{m,m}$, is an improper PGD having blocks of size 2.
The dual of a $K_{m,m}$, a quadratic grid of side $m$, is a PGD with parameters \[v=m^2,\ \ b=2m,\  \ k=m,\ \  r=2,\ \ \alpha =1,\ \  \beta = m+1.\] 
\end{example}
 
Observe that the PGD obtained from $K_{m,m}$ is TD$_1(2, m)$.
Transversal designs form a special family of group divisible designs in the following sense. 
 \begin{definition} \label{GDD}
 A tactical configuration, a $1$-$(v, b, k, r)$ design, is called a \emph{group divisible design} with parameters $(v, b, k, r, g, u; \mu_1, \mu_2)$, denoted by GDD$(v, b, k, r, g, u; \mu_1, \mu_2)$,  if the points can be divided into $g$ groups, each with $u$ points, so that $v=gu$ and two points belonging to the same group occur together in $\mu_1$ blocks and two points belonging to different groups occur together in $\mu_2$ blocks. 
  \end{definition}

The combinatorial properties of this family of GDDs have been studied in \cite{BC} and the analysis along with other types of designs is investigated in \cite{BSh} as early as in 1950s.
 The parameters of a GDD$(v, b, k, r, g, u; \mu_1, \mu_2)$ satisfy the basic relations 
\[v=gu, \ \ \mu_1(u-1)+\mu_2u(g-1) =r(k-1), \ \  r-\mu_1\ge 0, \ \ rk-v\mu_2\ge 0.\]

The combinatorial properties and some construction methods of the GDDs 
have been investigated by Bose, Shrikhande and Bhattacharya in \cite{BSB}. 

We note that a GDD$(v, b, k, r, g, u; \mu_1, \mu_2)$ with $\mu_1=\mu_2=\lambda$ is a $2$-$(v, k, \lambda)$ design. Also a GDD$(v, b, k, r, g, u; \mu_1, \mu_2)$ with $\mu_1=0$ reduces to a transversal design.
A particular example of such a design with maximal possible number of blocks without having any repeated blocks has been given by Olmez and Song in \cite{OS} as follows. 

\begin{example} Let $P$ be a $ql$-element set with a partition $\mathcal{P}$ of $P$ into $l$ parts
of size $q$. Let $\mathcal{P} = \{P_1, P_2, \dots, P_l\}$ and let
\[\mathcal{B} = \{B\subset P : |B\cap P_i| = 1,\ \  \forall i = 1, 2, \dots, l\}.\]
Then $\mathcal{B}$ consists of $q^l$ subsets of size $l$ as blocks. It is clear that the pair $(P,\mathcal{B})$ forms a tactical configuration with parameters $(v, b, k, r)=(ql, q^l, l, q^{l-1})$.
It has the property that any two points
from the same group never occur together in a block while any two points
from different groups occur together in $q^{l-2}$ blocks. It is a GDD$(ql, q^l, l, q^{l-1}, l, q; 0, q^{l-2})$ and a  
$TD_{\lambda}(l, q)$ with $\lambda=q^{l-2}$.
\end{example}

\subsection{Partial geometric difference sets}{\label{secPGDS}}
Olmez \cite{O} introduced the notion of a partial geometric difference set as a new source of a PGD. Then in \cite{NOS}, Nowak, Olmez and Song extended this notion to difference families.
\begin{definition}  Let $G$ be a group of order $v$. Let $S$ be a $k$-element subset of $G$ with $2 \leq k\leq  v$.
 For $g\in G$, we define $\delta_S(g)$ (or $\delta(g)$ in short if $S$ is given as fixed) by \[\delta(g):=|\{(s, t)\in S\times S: g=st^{-1}\}|. \]
 A $k$-subset $S$ of a group $G$ of order $v$ is called a \emph{partial geometric difference set} (PGDS) in $G$, denoted by PGDS$(v, k;\alpha, \beta)$ if there exist constants $\alpha$ and $\beta$ such that, 
 \[\forall x\in G, \quad \sum\limits_{y\in S}\delta(xy^{-1}) \ \ = \ \ \left \{\begin{array}{ll} \alpha & \mbox{if } \ x\notin S;\\ \beta & \mbox{if }\ x\in S.\end{array}\right .\]
 \end{definition}
 
\begin{definition} Let $v, k, n$ be positive integers with $v>k>2$. Let $G$ be a group of order $v$. 
Given a family  $\mathcal{S}=\{S_1, S_2, \dots, S_n\}$ of $n$ distinct non-$\emptyset$ $k$-element subsets of $G$, we define
\[\delta_{\mathcal{S}}(g):= \sum\limits_{j=1}^n \delta_{S_j}(g)=\sum\limits_{j=1}^n |\{(s, t)\in S_j\times S_j: g=st^{-1}\}|. \]
A family  $\mathcal{S}=\{S_1, S_2, \dots, S_n\}$ of $n$ distinct non-$\emptyset$ $k$-element subsets of $G$ is called  a \emph{partial geometric difference family} (PGDF) in $G$ with parameters $(v,k,n;\alpha, \beta)$ for some integers $\alpha$ and $\beta$, if 
\[\forall x\in G\ \ \mbox{and} \ \ \forall i\in\{1, 2, \dots, n\},\quad  \sum\limits_{r\in S_i}\delta_{\mathcal{S}}(xr^{-1})=\left\{\begin{array}{ll} \alpha & \mbox{if } x\notin S_i \\
\beta & \mbox{if } x\in S_i \end{array}\right .\]
\end{definition}

Olmez observed the following relationship between a PGDS and a PGD:
\begin{theorem}\cite{O} 
If $S$ is a PGDS$(v, k; \alpha, \beta)$ in a group $G$, then $(G, \mathcal{B})$ where $\mathcal{B}=\{gS: g\in G\}$ with block $gS=\{gs: s\in S\}$ is a PGD$(v, v, k, k; \alpha, \beta)$.
\end{theorem}
 
Known construction methods of PGDSs use group characters and group rings. Given a finite group $G$,  for any non-$\emptyset$ subset $S\subseteq G$, in the group ring 
$$\mathbb{Z}G =\left \{\sum\limits_{g\in G} z_g g:\  z_g\in\mathbb{Z}  \right \}\quad  (\mbox{the ring of formal sums of elements of $G$}),$$  
the element $\sum_{s\in S}s$ is denoted by  $\underline{S}$ and called \emph{simple quantity}; and so,  for $S^{-1}=\{s^{-1}: s\in S\}$, $\underline{S}^{-1}=\sum_{s\in S}s^{-1}$.  Olmez \cite{O} proved the following criteria for the existence of a PGDS:

 \begin{lemma} \cite{O} Let $S$ be a $k$-subset of a group $G$ of order $v$. Then 
 \begin{enumerate}
     \item $S$ is a PGDS$(v, k; \alpha, \beta)$ in $G$ if and only if the simple quantity $\underline{S}$ of $S$ satisfies $\underline{S} \ \underline{S}^{-1}\underline{S}= (\beta-\alpha) \underline{S} + \alpha\underline{G}$ in the group ring $\mathbb{Z}G$;
     \item $S$ is a PGDS$(v,k;\alpha, \beta)$, then $k(v-k)\equiv 0 \pmod{(\beta-\alpha)}$.
     \end{enumerate}
 \end{lemma}
 
\begin{theorem} \cite{NOS} Let $\mathcal{S}=\{S_1, S_2, \dots, S_n\}$ be a family of distinct $k$-subsets of a group $G$ of order $v$. If $\mathcal{S}$ is a PGDF$(v, k, n; \alpha, \beta)$, then $(G, \mathcal{B})$, where $\mathcal{B}=\bigcup\limits_{j=1}^n\{gS_j: g\in G\}$ is a 
PGD$(v, vn, k, kn; \alpha, \beta)$.
\end{theorem}
We now give several examples. Herein and in what follows, we denote the circulant matrix with the first row entries $[c_0, c_1, \dots, c_{v-1}]$ by $C[c_0, c_1, \dots, c_{v-1}]$.
\begin{example}\label{PGDS1}\cite{O}  Let $G=Q_8=\{1, i, j, k, -1, -i, -j, -k\}$, the quarternion group of order 8. It is shown that $S=\{-1, i, j, k\}$  is a PGDS in $G$. Then $(G, \mc{B})$ with $$\mc{B}= \{gS: g\in G\}=\{S, -S, iS, -iS, jS, -jS, kS, -kS\}$$ is a PGD$(8, 8, 4, 4; 6, 10)$.
The rows and columns of $N$ are indexed by the elements of $G$ and $\mc{B}$ in the order given in the corresponding sets, so that  we have concurrence
\[ \lambda_{xy} = \left \{\begin{array}{ll}  0 & \mbox{if } x=-y,\\ 2 &\mbox{if } x\neq \pm y;\end{array}\right .\qquad \mbox{and } \ 
NN^T= C[4, 2, 2, 2, 0, 2, 2, 2].\]
Notice that Spec$(NN^T) =[16^1, 4^4, 0^3]$. 
\end{example}

\begin{example}\label{PGDS2}  Let $G=A_4$, the alternating group of order 12. Let 
\[S_1 = \{(1), (234), (243), (12)(34), (123), (124)\},\ \ S_2=\{ (1), (234), (243), (12)(34), (132), (142)\}.\]
Then both $S_1$ and $S_2$ are PGDSs with the same parameters $(12, 6; 12, 24)$. 

The associated PGDs $\mathcal{D}_1=(A_4, \{gS_1:g\in A_4\})$ and $\mathcal{D}_2=(A_4, \{gS_2:g\in A_4\})$ have the same parameters 
$(12, 12, 6,  6; 12, 24)$ and the same spectrum Spec$(NN^T)=[36^1, 12^3, 0^8]$. However, they are not isomorphic as they have different concurrences. Namely, 
\[NN^T(\mathcal{D}_1) = C[6, 2, 2, 2, 6, 2, 2, 2, 6, 2, 2, 2]\ \mbox{ and }\ 
NN^T(\mathcal{D}_2)= C[6, 3, 3, 0, 3, 3, 6, 3, 3, 0, 3, 3].\]  
\end{example}

\begin{example}
\begin{enumerate}
\item In $\mathcal{D}_1$ in the above example, if we take only half of the blocks without containing repeated blocks with $S=S_1$, for example,
$\{S, g_1S, g_2S, g_4S, g_5S, g_6S\}$ as a block set with $G$ as point set, it is a PGD$(12, 6, 6, 3; 6, 12)$ with  the concurrence matrix $C=C[3, 1, 1, 1, 3, 1, 1, 1, 3, 1, 1, 1]$ and  Spec$(C)=[18^1, 6^3, 0^8]$.
Notice that it is a subdesign of $\mc{D}_1$ with the selected 6 blocks, the parameters except $v$ become halved from the original (but $\sigma$ remained as the same).
\item In $\mathcal{D}_2$, if we take selected four blocks (without containing repeated blocks) with $S=S_2$, for example,
$\{S, g_3S, g_6S, g_9S\}$ as a block set, we get a PGD$(12, 4, 6, 2; 4, 8)$ with its concurrence matrix  $C[2, 1, 1, 0, 1, 1, 2, 1, 1, 0, 1, 1]$ and its spectrum $[12^1, 4^3, 0^8]$.
\end{enumerate}
\end{example}

\begin{remark}
We have seen in the above example that a PGD may contain a `subdesign' (a design on the same point set but its block set is a subset of the given block set). This indicates a possibility of obtaining a larger PGD from a given PGD. This includes the method obtaining a PGD from a given PGD simply by repeating every block a constant number of times. That is, given a PGD $\mathcal{D}=(P, \mathcal{B})$ we have a new PGD $\mathcal{D}'=(P, \ l.\mathcal{B})$ by repeating each block $l$ times. In this case, we  will denote $\mathcal{D}'$ by $\mathcal{D}\otimes J_{1, l}$ as the incidence matrix $N'$ for $\mathcal{D}'$ can be expressed by  \[N'=N\otimes J_{1,l}=[\underbrace{N|N|\cdots |N}_l]\]
where $N$ is the incidence matrix of $\mathcal{D}$. There are various other ways to obtain new designs from old ones as we will see along the way.
\end{remark}

A new PGDS can be obtained from old one using the direct product of groups as below. This implies that there are infinitely many examples of PGDs coming from PGDSs.  
\begin{lemma} 
(\cite[Theorem 3.11]{O}) Given a PGDS$(v,k;\alpha, \beta)$ $S$ in a group $G$, the set $R=S\times \mbb{Z}_m\subset G\times\mbb{Z}_m$ is a PGDS$(vm, km; m^2\alpha, m^2\beta)$ in group $G\times\mbb{Z}_m$.
\end{lemma}
\begin{example}\cite{O} In group $\mbb{Z}_4$, the subset $S=\{0, 1\}$ is a PGDS$(4, 2; 1,3)$. By lemma,
the set $R:=S\times \mbb{Z}_3 \subset \mbb{Z}_4\times \mbb{Z}_3$ is a PGDS$(12, 6; 9, 27)$. From this PGDS, we obtain a PGD$(12, 12, 6, 6; 9, 27)$ with point set $$\mbb{Z}_4\times \mbb{Z}_3=\{(0,0), (0, 1), (0, 2), (1, 0), (1,1), (1,2), (2,0), (2,1), (2,2), (3,0), (3,1), (3,2)\}$$ and block set  $\mc{B}=\{(i,j)R: (i,j)\in \mbb{Z}_4\times \mbb{Z}_3\}$. Then the incidence matrix $N$ of PGD $(\mbb{Z}_4\times \mbb{Z}_3, \mc{B})$, and  $NN^T$ and $NN^TN$, respectively, are as the following:
\[\left [\begin{array}{cccc} J & 0 & 0 & J\\ J & J & 0 & 0\\ 0 & J & J & 0\\ 0 & 0 & J & J\end{array}\right ], \quad 
NN^T=\left [\begin{array}{cccc} 6J & 3J & 0 & 3J\\ 3J & 6J & 3J & 0\\ 0 & 3J & 6J & 3J\\ 3J & 0 & 3J & 6J\end{array}\right ], \quad
NN^TN= \left [\begin{array}{cccc} 27J & 9J & 9J & 27J\\ 27J & 27J & 9J & 9J\\ 9J & 27J & 27J & 9J\\ 9J & 9J & 27J & 27J\end{array}\right ],\]
where $J$ denotes $3\times 3$ all-ones matrix. We note that Spec$(NN^T)=[36^1, 18^2, 0^9]$. 
\end{example}
There are many new PGDSs and PGDFs reported lately (for example, see \cite{Mi, CCZ}).



\subsection{Strongly regular graphs} \label{sec-srgs}
Many researchers have reported interesting relationship between strongly regular graphs, symmetric 2-designs and partial geometric designs in various contexts. (For instance, see \cite{Br, CvL, GS, N, vDS}.) It is well-known that every complete multipartite strongly regular graph gives rise to a PGD. Recently, Nowak, Olmez and Song \cite{NOS1} have shown that a strongly regular graph with parameters $(\textbf{v},\textbf{k},\nu, \mu)$ gives rise to a PGD if and only if it satisfies either $\textbf{k}=\mu$ or $\nu=\mu$. They also have characterized which strongly regular graphs give rise to symmetric PGDs, where \emph{symmetric} means the number of blocks and the points are the same. We briefly recall what they did.

 A \textit{strongly regular graph} with parameters $(\textbf{v}, \textbf{k}, \nu, \mu)$ is a simple (undirected) graph $\Gamma$  that satisfies the conditions: (i) $\Gamma$ has $\textbf{v}$ vertices, (ii) $\Gamma$ is regular with degree (or valency) $\textbf{k}$, (iii) any two adjacent vertices have exactly $\nu$ common neighbors, and (iv) any two non-adjacent vertices have exactly $\mu$ common neighbors. Such a graph is denoted by SRG$(\textbf{v}, \textbf{k}, \nu,\mu)$.

A graph and its adjacency matrix will be used interchangeably.
We note that a $\{0,1\}$-matrix $A$ is the adjacency matrix of a SRG$(\textbf{v}, \textbf{k}, \nu,\mu)$ if and only if 
\begin{equation}\label{srg}
AJ=\textbf{k}J=JA, \quad A^2=\textbf{k}I+\nu A+\mu(J-I-A).\end{equation}

Herein,
the strongly regular graphs we use will be `nontrivial'; that is, neither a complete graph nor a discrete graph.
We will say that ``a \emph{graph $\Gamma=(X, R)$ with its vertex set $X$ and edge set $R$ gives rise to a design $(P, \mathcal{B})$}" if the adjacency matrix $A$ of $\Gamma$ is \textit{equivalent} to the incidence matrix $N$ of $(P, \mathcal{B})$.
That is, for each $v\in X$, if we let $N_v=\{x\in X\colon \{x, v\}\in R\}$ and $\mathcal{N}=\{N_v\colon v\in X\}$, the pair $(X, \mathcal{N})$ forms a design that is isomorphic to $(P, \mathcal{B})$. 
\begin{lemma} 
\cite[Lemma 2]{NOS1} \label{srg-dsgn}
Let $\Gamma$ be a $\rm{SRG}(\textbf{v},\textbf{k},\nu,\mu)$. Let $A$ be the adjacency matrix of $\Gamma$. Then $A^3=\beta A+\alpha (J-A)$ for some integers $\alpha$ and $\beta$ if and only if either $\nu=\mu$ or $\textbf{k}=\mu$. 
\end{lemma}

Every complete multipartite strongly regular graph can be viewed as the complement of $c$ copies of the complete graph $K_m$ on $m$ vertices for some integers $c$ and $m$ (where $c, m\ge 2$).  We denote such a graph by $\overline{cK_m}$.

\begin{cor} \cite[Corollary 3]{NOS1}
The complete multipartite strongly regular graph $\overline{cK_m}$ gives rise to a PGD with parameters \[(cm, \ \ cm,\ \  (c-1)m,\ \ (c-1)m;\ \  (c^2-3c+2)m^2,\ \ (c^2-3c+3)m^2).\]
\end{cor}

A $\mbox{SRG}(\textbf{v},\textbf{k},\nu, \nu)$ is also called a $(\textbf{v},\textbf{k},\nu)$-graph. The adjacency matrix $A$ of a $(\textbf{v},\textbf{k},\nu)$-graph satisfies the identity $A^2=\textbf{k}I+\nu (J-I)$; therefore, it gives rise to a symmetric 2-$(\textbf{v},\textbf{k},\nu)$ design. 

\begin{cor}\cite[Corollary 4]{NOS1}
A $(\textbf{v},\textbf{k},\nu)$-graph gives rise to a PGD with parameters \[(\textbf{v}, \textbf{v},\textbf{k}, \textbf{k};\textbf{k}\nu, \textbf{k}\nu+\textbf{k}-\nu).\]
\end{cor}

\begin{example} \label{rem-srg}
It is well-known that both the Hamming graph $H(2,4)$ and the Shrikhande graph are $(16, 6, 2)$-graphs. Although they are non-isomorphic strongly regular graphs,  they give rise to the same 2-$(16, 6, 2)$-design.  Hence, we have a PGD with parameters $(\textbf{v}, \textbf{v},\textbf{k}, \textbf{k};\alpha, \beta)=(16, 16, 6, 6; 12, 16)$.
\end{example}


\subsection{Association schemes}
It is known by Nowak, Olmez and Song \cite{NOS1} that PGDs are also coming from a certain family of association schemes of class 3. 
They investigated the characteristics of the graphs that give rise to PGDs and observed that some of these graphs arise as the relation graphs of 3-class association schemes. A family of such 3-class association schemes are the 3-class `fusion schemes' of $d$-class Hamming schemes over the finite field of order 3 for $d\ge 3$. 

\begin{definition} Let $X$ be a finite set and let $\{R_i\subseteq X\times X: i\in [d]\}$, where $[d]:=\{0,1,\ldots, d\}$, be the set of relations on $X$. Then the pair $(X, \{R_i
\}_{0\le i\le d})$ is called an association scheme
of class $d$ if the following conditions are satisfied. \medskip
\begin{enumerate}
\item[(i)]  $R_0=\{ (x,x):\ x\in X\} $, 
\item[(ii)] $ R_0 \cup R_1
\cup \cdots \cup R_d = X\times X$ and $R_i\cap R_j = \emptyset $ if
$i\neq j$,
\item[(iii)] For each $i\in [d] $,
if we define ${}^tR_i=\{ (y,x) :\ (x,y) \in R_i\} $, then
${}^tR_i=R_i$.\footnote{Symmetric associations are considered in this paper.}
\item[(iv)] For $h, i,j\in [d] $, the number $|
\{ z\in X :\ (x,z)\in R_i, (z,y)\in R_j\}| $ does not depend on the
choice of $x,y$ under the condition $ (x,y)\in R_h$, and is equal to
$p_{ij}^h$ which depends only on $h, i, j$.
\end{enumerate}
\end{definition}
The integers $p_{ij}^h$, $h, i, j\in [d]$, are called the \textit{intersection numbers} or parameters of $\mathcal{X}$. These parameters are conveniently presented by $d+1$ matrices, called the intersection matrices $B_0, B_1, \dots, B_d$, where the $(j, h)$-entry of $B_i$ is:
\[\left [B_i\right ]_{jh} :=  p_{ij}^h.\]

Given a $d$-class association scheme $(X, \{R_i\}_{0\le i\le d})$, let $A_i$ be the adjacency matrix
with respect to the relation $R_i$. That is, we define $(A_i)_{x,y}=
1$, if $(x,y)\in R_i$; $(A_i)_{x,y}=0$, if $(x,y)\notin R_i.$ Then by the definition of the association scheme we have 
\begin{enumerate}
\item[(a)] $A_0=I$, $A_0+A_1+\cdots +A_d=J$, where $J$ is the all-ones matrix,
\item[(b)] $A_i=A_i^T$ for each $i\in[d]$, and
\item[(c)] for any $h, i, j\in [d]$, there exists a constant $p_{ij}^h$ such that
\begin{equation} A_iA_j=\sum\limits_{h=0}^d p_{ij}^h A_h.\end{equation}
\end{enumerate}
Therefore, the vector space spanned by $A_0,A_1,\ldots ,A_d$ over $\mathbb{C}$ is of
dimension $d+1$ and is closed under the ordinary matrix multiplication. 
That is, the adjacency matrices span an algebra $\mathcal{A}=\langle A_0, A_1, \dots, A_d\rangle$ of symmetric matrices over the complex, and it is called the \emph{Bose-Mesner algebra} of $\mathcal{X}$. 
As $\mathcal{A}$ is known to be a semisimple algebra, it admits another basis consisting of the central primitive idempotents:  $\{E_j\}_{j\in [d]}$ with $E_0=\frac{1}{n}J$. 

Let $P=\left [p_j(i)\right ]$ and $Q=\left [q_j(i)\right ]$ be the $(d+1)\times (d+1)$ base-change
matrices of $\mathcal{A}$ whose $(i,j)$-entries $p_j(i)$ and $q_j(i)$ are defined by
\begin{equation}
A_j=\sum\limits_{i\in [d]} p_j(i)E_i,\qquad
E_j=\frac{1}{n}\sum\limits_{i\in [d]}q_j(i)A_i,\quad j\in [d].
\end{equation}

We call $P$ the \emph{character table} of $\mathcal{X}$. Clearly $PQ=QP=nI$.
We note that
\begin{equation} \label{eig}
A_i E_j = p_i(j) E_j \hspace{1in}  (i,j \in [d]).
\end{equation}
In particular,  $p_i (j)$ is the eigenvalue of $A_i$ associated with the eigenspace $E_j \mathbb{C}^X$
$(\forall i,j\in [d])$.

We are now ready to describe the parameter sets of certain 3-class association schemes that give rise to PGDs. In \cite{NOS1}, it is shown that if a 3-class association scheme of order $3m^2$ satisfies certain parametric conditions, then its adjacency matrices $A_0, A_1, A_2, A_3$ satisfy the following identities for some constants $\alpha_i$ and $\beta_i$:
\begin{equation}\begin{array}{c}
A_1^3=\beta_1A_1 +\alpha_1(J-A_1),\\
A_2^3=\beta_2A_2+\alpha_2(J-A_2),\\
(A_3+A_0)^3=\beta_3(A_3+A_0)+\alpha_3(J-A_3-A_0).\end{array}\end{equation}

We state the characterization theorem of 3-class association schemes whose relation graphs give rise to PGDs.
\begin{theorem}(\cite{NOS1}) \label{thm-z} Let $\mathcal{Z}$ be a 3-class association scheme, and let $A_0, A_1, A_2, A_3$ be its adjacency matrices. Suppose that the character table $P$ of $\mathcal{Z}$ is given by
\[P =\left [ \begin{array}{cccc}
1 & m(m-1) & m(m+1)& (m-1)(m+1)\\
1 &  m   &  0&         -m-1  \\
1 & 0   & -m  &  m-1 \\
1 & -m&  m & -1      \\
 \end{array}\right ]. \]
Then the relation graphs $A_1, A_2, \mbox{ and } A_3+A_0$ of $\mathcal{Z}$ give rise to three symmetric ($v=b$) PGDs. In this case,
the parameters $(v, k; \alpha, \beta)$ of corresponding PGDs are given by
\[(3m^2,\ m(m-1);\ \frac13 m^2(m^2-3m+2),\ \frac13 m^2(m^2-3m+5)),\]
\[(3m^2,\ m(m+1);\ \frac13 m^2(m^2+3m+2),\ \frac13 m^2(m^2+3m+5)),\]
\[(3m^2,\ m^2;\ \frac13 m^2(m^2-1),\ \frac13 m^2(m^2+2)).\]
\end{theorem}

We recall the Hamming scheme $H(d,q)$. 
\begin{definition} Let $F$ be a $q$-element set and let \[X:=F^d=\{(x_1, x_2,\dots, x_d): x_j\in F, \ j=1, 2, \dots, d\}.\]
Define the association relation between any $\textbf{x}=(x_1, x_2, \dots, x_d)$ and $\textbf{y}=(y_1, y_2,\dots, y_d)$ in $X$ according to the Hamming distance \[\delta(\textbf{x}, \textbf{y}):=|\{j\in \{1, 2,\dots, d\}: x_j\neq y_j\}|;\] that is, define \[(\textbf{x},\textbf{y})\in R_i\ \Leftrightarrow \ \delta(\textbf{x}, \textbf{y})=i.\]
Then
$(X, \{R_i\}_{0\le i\le d})$ is an association scheme called the $d$-class \emph{Hamming scheme} over $F$, denoted by $H(d, q)$.\end{definition}

\begin{example} The association scheme whose eigenmatrices $P, Q$ are given below is known as 
the Hamming association scheme, $H(3,3)$, with $B=B_1$ ($B_{jh}=p_{1j}^h$).
\[P=Q =\left [ \begin{array}{cccc}
1 & 6 & 12 & 8\\
1 & 3 & 0 & -4\\
1 & 0 & -3 & 2\\
1 & -3 & 3 & -1\\
\end{array}\right ] \qquad
B=\left [ \begin{array}{cccc}
0 & 1 & 0 & 0\\
6 & 1 & 2 & 0\\
0 & 4 & 2 & 3 \\
0 & 0 & 2 & 3\\ \end{array}\right ] \]
Through direct calculation,  we have
\[A_1^3=15A_1+6(J-A_1),\quad 
 A_2^3=69A_2+60(J-A_2),\] \[ (A_3+I)^3=33(A_3+I)+24(J-A_3-I).\]
\end{example}

The infinite family of PGDs are now obtained as the following.
\begin{theorem} \label{sch-f}\cite{NOS1} Let $(X, \{R_i\}_{0\le i\le d})$ be the Hamming scheme $H(d,3)$ with $d=2l+1$ for  $l\ge 1$. Let $S_0=R_0$ and define
\[S_j=\bigcup\limits_{i=0}^{\lfloor (d-j)/3\rfloor }R_{3i+j}, \mbox{ for } j=1, 2, 3\]
where $\lfloor (d-j)/3\rfloor$ denotes the greatest integer less than or equal to $(d-j)/3$.
Then $\mathcal{F}=(X, \{S_0, S_1, S_2, S_3\})$ is a 3-class scheme.
 For each $l\ge 1$, the relation graphs of the 3-class scheme $\mathcal{F}$ above give rise to three non-isomorphic symmetric PGDs with parameters $(v, k;\alpha, \beta)$:
\[\big(3^{2l+1},\ 3^{2l}+(-1)^l3^l;\ 3^{4l-1}+(-1)^l3^{3l}+2\cdot 3^{2l-1}, \ 3^{4l-1}+(-1)^l3^{3l}+5\cdot 3^{2l-1}\big),\]
\[\big(3^{2l+1},\  3^{2l}-(-1)^{l}3^l; \ 3^{4l-1}-(-1)^{l}3^{3l}+2\cdot 3^{2l-1},\ 3^{4l-1}-(-1)^{l}3^{3l}+5\cdot 3^{2l-1}\big),\]
\[\big(3^{2l+1}, \ 3^{2l}; \ 3^{4l-1}-3^{2l-1}, \ 3^{4l-1}+2\cdot 3^{2l-1}\big).\]
\end{theorem}


\subsection{Partial geometries}

Partial geometries were introduced by Bose \cite{Bose}  in order to provide a generalization for known characterization theorems for strongly regular graphs. They are important sources of PGDs. There are several surveys on partial geometries: e.g., see Cameron and van Lint \cite {CvL} and references in there. 
A finite incidence structure is called a \emph{partial plane} (or \emph{geometric}) if it satisfies one of the two equivalent conditions: Two distinct points have at most one common incident block, and two distinct blocks have at most one common incident point. 
So we recall that a partial geometry $PG(\kappa, \rho,\tau)$ is a geometric tactical configuration with parameters $(v, b, \kappa, \rho)$ such that for any antiflag $(a, B)$, there are exactly $\tau\ge 1$ lines containing point $a$ and intersecting block $B$.

\begin{definition}
A \textit{partial geometry} $PG\left(\kappa, \rho, \tau\right)$ is a set of points $P$, a set of lines $\mathcal{L}$, and an incidence relation between points and lines satisfying the following properties:
\begin{enumerate}
    \item Every line is incident with exactly $\kappa$ points $\left(\kappa \geq 2\right)$, and every point is incident with exactly $\rho$ lines $\left(\rho \geq 2\right)$.
    \item Any two points are incident with at most one line.
    \item If a point $p$ and a line $L$ are not incident, then there exists exactly $\tau$ lines $\left(\tau \geq 1\right)$ that are incident with $p$ and meet $L$.
\end{enumerate}\end{definition}


The following properties and facts on partial geometries are found in \cite{BSS}:
\begin{prop} 
\begin{enumerate}
\item A $PG(\kappa, \rho, \tau)$ is a tactical configuration with parameters \[\left(v, b, k, r\right)  = \left(\kappa+\frac{1}{\tau}\kappa(\kappa-1)(\rho-1),\ \  \rho+\frac{1}{\tau}\rho(\rho-1)(\kappa-1),\ \ \kappa, \ \ \rho \right).\] 
\item A $PG(\kappa, \rho, \tau)$ is a 2-$(v, \kappa, 1)$ design if and only if $\tau=\kappa$.\\
A $PG(\kappa, \rho, \tau)$ is the dual of  a 2-$(v, \kappa, 1)$ design if and only if $\tau=\rho$.\\
A $PG(\kappa, \rho, \tau)$ is a transversal design with $\lambda=1$ if and only if $\tau=\kappa-1$.\\
A $PG(\kappa, \rho, \tau)$ is a net if and only if $\tau=\kappa$.

\item For a $PG(\kappa, \rho, \tau)$, the graph $\Gamma=(X, R)$ given by $xy\in R$ iff $x$ and $y$ are incident with the same line, is a strongly regular graph, and some eigenvalue of $\Gamma$ has multiplicity
\[\frac{\kappa(\kappa-1)\rho(\rho-1)}{\tau(\kappa+\rho-\tau-1)}.\]
\end{enumerate}
\end{prop}

\begin{lemma} Let $\mathcal{P}$ be a tactical configuration with parameters $(v, b, k, r)$.
\begin{enumerate}
\item[(a)] $\mathcal{P}$ is a partial plane iff $s(a, B) = r+k-1$ for every flag $(a, B)$. 
\item[(b)] A partial plane is a PGD iff it is a partial geometry.
\item[(c)] For a PGD $\mathcal{P}$, the following three are equivalent: 

  \begin{tabular}{l}
 (i) $\beta=r+k-1$,\\
 (ii) $\mathcal{P}$ is partial plane,\\
 (iii) $\mathcal{P}$ is a partial geometry.
  \end{tabular}
\end{enumerate}
\end{lemma}

\begin{prop} A partial geometry with parameters $(\kappa, \rho, \tau)$ is a PGD with parameters
\[\left (\frac{\kappa (\tau+(\kappa-1)(\rho-1))}{\tau}, \ \frac{\rho (\tau+(\kappa-1)(\rho-1))}{\tau}, \  \kappa, \ \rho; \ \ \tau, \  \kappa+\rho-1\right )\]
and its concurrence satisfies \[\lambda_{xy} = \left \{\begin{array}{ll} 0 & \mbox{if } x \mbox{ and } y \mbox{ not collinear}; \\ 1 & \mbox{if } x \mbox{ and } y \mbox{ are collinear}.\end{array}\right .\]
\end{prop}

\begin{example}\label{PG}
Consider the incidence structure defined on the set $\mathbb{Z}_q\times \mathbb{Z}_q$ (where $q$ is prime) as the following:
\begin{itemize}
\item its point set: $P=\{p=p(x,y): x,y\in\mathbb{Z}_q\}$;
\item its blocks (lines): $\mathcal{B}=\{l=l(k,j): y=kx+j\}$;
\item incidence: point $p(x,y)$ is incident to line $l(k,j)$ if and only if $y=kx+j$.
\end{itemize}

Suppose $n=3$. Then
\[P=\{ p_1(0,0),  p_2(1,0), p_3(2, 0), p_4(0,1), p_5(1, 1), p_6(2,1),
 p_7(0,2), p_8(1, 2), p_9(2,2)\};\]
\[\mathcal{B} = \left \{\begin{array}{ccc}
l_1(0,0)=\{p_1, p_2, p_3\}, &  l_2(0, 1)=\{p_4, p_5, p_6\}, &  l_3(0,2)=\{p_7, p_8, p_9\},\\
l_4(1,0)=\{p_1, p_5, p_9\}, &  l_5(1, 1)=\{p_4, p_8, p_3\}, &  l_6(1, 2)=\{p_7, p_2, p_6\},\\
l_7(2,0)=\{p_1, p_8, p_6\}, &  l_8(2,1)=\{p_4, p_2, p_9 \}, &  l_9(2, 2)=\{p_7, p_5,p_3\}.\end{array} \right \}\]

We observe the following:
\begin{enumerate}
\item $(P, \mathcal{B})$ is a partial geometry PG$(3,2,2)$, and a PGD$(9,9,3,3;2,5)$. 
\item The concurrence matrix $NN^T$ is the circulant matrix $C[3,1,1,0,1,1,0,1,1]$. 
It is also shown that the spectrum of $NN^T$ is $[9^1, 3^6, 0^2]$. 
\end{enumerate}
\end{example}

\subsection{Special partially balanced incomplete block designs}{\label{SPBIBD}}
Both partially balanced and special partially balanced incomplete block designs have been sought in the construction of PGDs from the time  when PGDs were introduced.
A \emph{partially balanced incomplete block design} (PBIBD) may be defined by replacing the balanced condition in the notion of BIBD by a weaker condition. That is, in a PBIBD, every pair of points do not necessarily appear together equally often, but occur together in a certain number of blocks depending on an underlying association relation between the points (see Raghavarao \cite{Ra}). So given a $v$-element set $X$, the design is defined based on the underlying association scheme on $X$ as in the following:
\begin{definition}
Let $\mathcal{X}=(X, \{R_i\}_{0\le i\le d})$ be a $d$-class association scheme defined on $X$. The pair $(X, \mathcal{B})$ where $\mathcal{B}$ is a collection of $b$ subsets of size $k$, is a PBIBD with $d$ associate classes attached to $\mathcal{X}$ if it satisfies:
\begin{enumerate}
    \item every point occurs at most once in a block,
    \item each point appears in exactly $r$ blocks, and
    \item any two points that are $i$th associates occur together in $\lambda_i$ blocks. 
\end{enumerate}
The numbers $v, b, k, r; \lambda_1, \lambda_2, \dots, \lambda_d; k_1, k_2, \dots, k_d$, where $k_i$ is the number of $i$th associates of a point,  are called the parameters of the PBIBD.
\end{definition}
We note that for a given association scheme, many different PBIBDs can be attached to the association scheme. It has been shown by Bose and Shimamoto in \cite{BSh}, all PBIBDs with two associate classes can be divided into a small number of types according to the nature of the association relations. One simple and important type is the GDD in Definition \ref{GDD}. 

\begin{definition} A \emph{special partially balanced incomplete block design of type $(s,t)$} is a tactical configuration $\mathcal{D}=(P, \mathcal{B})$ with parameters $(v, b, k, r)$ such that 
\begin{enumerate}
\item for some integers $\lambda_1>\lambda_2\ge 0$, every distinct pair of points occur together in either $\lambda_1$ or $\lambda_2$ blocks, 
\item for any point-block pair $(a, B)\in P\times \mathcal{B}$,
 \[|\{b\in B: a\neq b, \lambda_{ab}=\lambda_1\}|=\left \{\begin{array}{ll} s & \mbox{if } a\in B\\ t & \mbox{if } a\notin B\end{array}\right .\]
 where $s,t$ are nonnegative integers.
 \end{enumerate}
 \end{definition}
 
 Throughout, we denote the above design $\mathcal{D}=(P, \mathcal{B})$ by SPBIBD$(v, b, k, r; \lambda_1, \lambda_2)$ of type $(s, t)$.
The following facts have been discussed in Bridges and Shrikhande \cite{BSh}.
 \begin{remark}
 \begin{enumerate}
 \item An SPBIBD of type $(s,t)$ is a PBIBD with two associate classes with the additional property that for any point-block pair $(a, B)$ the number of first associates of $a$ in the block $B$ is $s$ or $t$ depending on whether $a\in B$ or $a\notin B$, respectively (cf. \cite{Bo} for PBIBD).
 \item A $2$-$(v, k, \lambda)$ design, is an SPBIBD  where $\lambda_1=\lambda$, 
$\lambda_2=0$, $s=k-1$ and $t=k$.
\item A PG$(\kappa, \rho, \tau)$ is an SPBIBD of type $(\kappa-1, \tau)$. It is known that the SPBIBDs of type $(k-1, t)$ with $\lambda_1=1$ and $\lambda_2=0$ are precisely partial geometries.  Further connections to various type of PBIBDs have been established by Bose and Conner, Raghvarao, and Bridges and Shrikhande (see \cite{Ra, BSh} and references in there).  
\end{enumerate}
\end{remark}

There are other examples of PGDs coming from SPBIBDs that are not from any other sources listed above. Some of them are given by van Dam and Spence in \cite{vDS} including the following. The concurrence matrices of these PGDs hardly appear as circulant.
\begin{remark} 
\begin{enumerate}
    \item \cite[Proposition 4]{vDS} \quad Let $\mathcal{D}$ be a PGD$(v, b, k, r; \alpha, \beta)$. Suppose Spec$(NN^T) =[kr^1, n^{v-2}, 0^1]$. Then $\mathcal{D}$ is an SPBIBD$(v, b, k, r; \lambda_1, \lambda_2)$. Furthermore, this SPBIBD is based on the 2-class association scheme one of whose relation graphs is $K_{v/2,v/2}$ where two points in the same part of the bipartition are contained in $\lambda_2=\frac{\alpha}{k}-\frac{n}{v}$ blocks and two points in different parts are contained in $\lambda_1=\frac{\alpha}{k}+\frac{n}{v}$ blocks.
    
    Examples of such PGDs are PGD$(9, 21, 3, 7; 5, 11)$ and PGD$(12, 20, 6, 10; 27, 33)$. 
    \item \cite[Proposition 7]{vDS} A PGD$(v, b, k, r;\alpha, \beta)$ with $k=3$ is an SPBIBD with $\lambda_1=\frac12(\beta-k-r+1)+1$
and $\lambda_2=0$ or a $2$-$(v, 3, \frac12(\beta-k-r+1)+1)$ design. Moreover, in the case of an SPBIBD, $\frac12(\beta-k-r+1)+1$ divides both $r$ and $n-r$.

Examples of such PGDs in this case are PGD$(9,9,3,3; 1, 7)$, PGD$(9, 15, 3, 5; 3, 9)$, PGD$(9, 18, 3, 6; 2, 14)$, and PGD$(9, 21, 3, 7; 5, 11)$.

\end{enumerate}
\end{remark}
 
 
 
 \subsection{Affine resolvable designs}

When we replace each point of a PGD$(v_0, b_0, k_0, r_0; \alpha_0, \beta_0)$ $\mathcal{D}$ with a set of $l$ points, we obtain a new PGD$(v=lv_0, b=b_0, k=lk_0, r=r_0;\alpha= l\alpha_0, \beta=l\beta_0)$, denoted by $\mathcal{D}\otimes J_{l,1}$. (As we have seen earlier, by repeating blocks $m$ times we obtain a new PGD, $\mathcal{D}\otimes J_{1,m}$, with parameters $(v_0, mb_0, k_0, mr_0; m\alpha, m\beta)$.)

Let $\mathcal{D}=(P,\mathcal{B})$ be a 2-$(v, k, \lambda)$ design. 
A \emph{parallel class} in $\mathcal{D}$ is a set of pairwise disjoint blocks that partition $P$; and so, $k|v$ and $r|b$ in this case. 

 \begin{definition} A 2-$(v, k, \lambda)$ design is said to be \emph{resolvable} if the block set $\mathcal{B}$ can be partitioned into $r$ disjoint parallel classes. A resolvable design is called \emph{affine} if there exists a positive integer $m$ such that any two non-parallel blocks intersect in exactly $m$ point(s). 
 \end{definition}
 
 All parameters of an affine (resolvable) design may be expressed in terms of $m$ and another parameter $q$ where $q=\frac{v}{k} = \frac{b}{r}$, the number of blocks in a parallel class, as follows:
 \[v=mq^2, \ b=q\frac{(mq^2-1)}{q-1},\ r = \frac{mq^2-1}{q-1},\ k=mq, \ \lambda = \frac{mq-1}{q-1}.\]
 
 Let $\mathcal{D}=(P,\mathcal{B})$ be an affine resolvable 2-$(v, b, k, r, \lambda)$ design with $m=\frac{k^2}{v}$ and $q=\frac{v}{k}$. We denote this design by $\mathcal{AD}_m(q)$. Let $\mathcal{AD}_m(q)^{(l)}$ denote the incidence structure obtained from $\mathcal{AD}_m(q)$ by considering all points in $P$ and taking the blocks of any $l$ parallel classes of $\mathcal{D}$ for given $1< l<r$. Let $\mathcal{B}^{(l)}$ denote the set of blocks in the selected $l$ parallel classes. Then $\mathcal{AD}_m(q)^{(l)}=(P, \mathcal{B}^{(l)})$ inherits the following properties from $\mathcal{AD}_m(q)$:
 \begin{enumerate}
 \item every point is incident with $l$ blocks, and every block is incident with $k$ points,
 \item any two points are incident with at most $\lambda$ blocks, and
 \item for each non-incident point-block pair $(p, B)$, there are exactly $l-1$ blocks containing $p$ each of which intersect $B$ in $m$ points.
 \end{enumerate}
 
 \begin{remark} If $m=1$, the $\mathcal{AD}_1(q)$ is the affine plane of order $q$. In this case,  $\mathcal{AD}_1(q)^{(l)}=(P, \mathcal{B}^{(l)})$ is a partial geometry $PG(q, l, l-1)$.  In particular, if $l=q$, it is a $PG(q, q, q-1)$ which is a transversal design $TD_1(q, q)$. Each group consists of $q$ points that were collinear in the affine plane, but not in  $\mathcal{AD}_1(q)^{(q)}$. Note that every pair of points in $P$ is contained in either exactly one group or in exactly one block, but not both. 
 \end{remark}
 
 
 \subsection{Classical finite geometries}
 
 We introduce the construction of a PGD using totally isotropic subspaces of the symplectic geometry over $\mathbb{F}_q$ described by Feng and Zeng in \cite{FZ}. Similar constructions of PGDs were given in \cite{FZZ} using the totally isotropic subspaces of unitary and orthogonal geometries.
 
 Let $K$ be a $2l\times 2l$ nonsingular alternate matrix. The set of all $2l\times 2l$ matrices $T$ over $\mathbb{F}_q$ satisfying $TKT^T=K$ forms the symplectic group $Sp_{2l}(q)$ of degree $2l$ over $\mathbb{F}_q$. 
 We denote the $2l$-dimensional vector space over $\mathbb{F}_q$ by $\mathbb{F}_q^{(2l)}$. By symplectic space we mean the vector space $\mathbb{F}_q^{(2l)}$ together with the action $\mathbb{F}_q^{(2l)}\times Sp_{2l}(q)\rightarrow \mathbb{F}_q^{(2l)}$ of $Sp_{2l}(q)$ on row vectors by $(\textbf{x}, T)\mapsto \textbf{x}T$.
 
 We will consider any $m$-dimensional subspace $U$ of $\mathbb{F}_q^{(2l)}$ as an $m\times 2l$ matrix $U$ whose rows are the vectors of a basis of $U$. Two matrices $U_1$ and $U_2$ of rank $m$ represent the same subspaces if and only if there exists an $m\times m$ nonsingular matrix $Q$ such that $U_1=QU_2$. A subspace $U$ is said to be of type $(m,s)$ if $\dim(U)=m$ and the rank of $UKU^T$ is $2s$. It is well-known that subspaces of type $(m,s)$ exist in $\mathbb{F}_q^{(2l)}$ if and only if $2s\le m\le l+s$. The subspaces of type $(m,0)$ are called the $m$-dimensional totally isotropic subspaces. It is clear that all 1-dimensional subspaces are totally isotropic. 
 
 \begin{theorem}(\cite{FZZ}) Let $P=S(n, 2l; q)$ and $\mathcal{B}=S(m, 2l; q)$  be the set of all $n$-dimensional and $m$-dimensional, respectively, totally isotropic subspaces in the $2l$-dimensional symplectic space  $Sp_{2l}(q)$ where respectively $1\le n<m\le l$. For $x\in P$ and $B\in \mathcal{B}$, define $x\in B$ if and only if $x\subset B$ as subspaces of $\mathbb{F}_q^{(2l)}$. Then the incidence structure $\mathcal{D}=(P, \mathbb{B}, \in )$ is a PGD if and only if $n=1$ and $m=l$. Furthermore, for $l\ge 3$,  the parameters  of $\mathcal{D}$ 
 are 

 \[v=\frac{q^{2l}-1}{q-1}, \ \ b=\prod\limits_{i=1}^l (q^i+1), \ \ k=\frac{q^{l}-1}{q-1}, \ \ r=\prod\limits_{i=1}^{l-1} (q^i+1),\]
 \[\alpha=\frac{q^{l-1}-1}{q-1}\prod\limits_{i=1}^{l-2} (q^i+1),   \ \ \beta=\frac{q(q^{l-1}-1)}{q-1}\left (\prod\limits_{i=1}^{l-2} (q^i+1)-1\right).\]
 
If $l=2$, then it is a symmetric PGD with $(v, k;\alpha, \beta)=((q+1)(q^2+1), q+1; 1, 2q+1)$.
 \end{theorem}

\subsection{Construction of PGDs with prescribed automorphisms} 
In \cite{NO}, Nowak and Olmez extended the method of constructing ordinary designs using group actions introduced by Kramer and Mesner \cite{KM} to construct PGDs. We describe their construction in this subsection together with a few examples.

Given integers $v$ and $k$ with $v>k>2$, let $X$ be a $v$-element set, and let $G\le S_X$ be a permutation group on $X$. Consider the naturally extended action of $G$ on $\mathcal{P}_k(X)={X\choose k}$, the collection of all $k$-subsets of $X$. Let $\mathcal{O}_1, \mathcal{O}_2, \dots, \mathcal{O}_s$ be the orbits of $G$ on $\mathcal{P}_k(X)$. Let $\iota= (i_1, i_2, \dots, i_q)$ be an increasing subsequence of subindices of the orbits; i.e., $1\le i_1<i_2<\cdots < i_q\le s$. For each $h$, $1\le h\le q$, and any $k$-set $Y_h\in \mathcal{O}_{i_h}$, define the $q\times v$ matrix $M_h$ whose $(u,t)$-entry is given by
\[(M_h)_{ut} = \sum\limits_{A:\ x_t\in A\in \mathcal{O}_{i_u}}|A\cap Y_h|.\]
Then $M_h$ is independent from the choice of the orbit representative $Y_h$ for $\mathcal{O}_{i_h}$ and uniquely determined up to the ordering of elements of $X$ (a permutation of columns). That is, if $M_h$ and $M_h^{\prime}$ are the matrices corresponding to two $k$-sets $Y$ and $Y'$ belonging to $\mathcal{O}_{i_h}$, respectively, then $M_h$ and $M_h^{\prime}$ are the same up to a permutation of the columns. Furthermore, it was observed by Nowak and Olmez that

\begin{lemma} 
Given an $\iota=(i_1, i_2, \dots, i_q)$ and $Y, Y'\in\mathcal{O}_{i_h}$ as above, for any $q$-vector $\mathbf{w}\in \mathbb{N}^q$,
\[\left (\mathbf{w}M_h\right )_t = \left \{\begin{array}{ll} \alpha & \mbox{if } x_t\notin Y\\ \beta & \mbox{if } x_t\in Y \end{array}\right . 
\quad \Leftrightarrow \quad \left (\mathbf{w}M_h^{\prime}\right )_t = \left \{\begin{array}{ll} \alpha & \mbox{if } x_t\notin Y'\\ \beta & \mbox{if } x_t\in Y' \end{array}\right . .\]
\end{lemma}

Then they gave the following existence criterion for PGDs by generalizing the well-known Kramer-Mesner Theorem for 2-designs.
\begin{theorem} (\cite[Theorem 6]{NO}) 
Let $G$ be a subgroup of the symmetric group $S_v$. Let $\mathcal{O}_1, \mathcal{O}_2, \dots, \mathcal{O}_s$ be the orbits of $G$ on $\mathcal{P}_k(X)$.
Let $L$ be the $s\times v$ matrix whose $(i,t)$-entry $(L)_{it}$ is given by
\[(L)_{it}=|\{A\in\mathcal{O}_i: x_t\in A\}|.\]  
There exists a PGD$(v, b, k, r;\alpha, \beta)$ having $G$ as an automorphism group
if and only if \\
(i) there exists a nonnegative integral solution $\mathbf{z}=[z_1, z_2, \dots, z_s]$, $z_i\in\mathbb{N}\cup \{0\}$, to the equation
\[\mathbf{z}L = r\mathbf{1}_v, \quad (\mathbf{1}_v: \mbox{ all-ones vector}),\] and (ii) if $\iota=(i_1, i_2, \dots, i_q)$ is the sequence of coordinate places $i$ where $z_i>0$ and  $\mathbf{w}=[z_{i_1}, z_{i_2}, \dots, z_{i_q}]$, then for all $h$ with $1\le h\le q$,
\[\left (\mathbf{w}M_h\right )_t = \left \{\begin{array}{ll} \alpha & \mbox{if } x_t\notin Y_h,\\ \beta & \mbox{if } x_t\in Y_h. \end{array}\right . \]
\end{theorem}

We now have a few particular examples that illustrate the use of the above criterion.
\begin{enumerate}
\item[(a)] (\cite[Ex.7]{NO})  Let $X=\{1, 2, 3, 4, 5, 6\}$. When the group $\mathbb{Z}_6$ acts on the set ${X\choose 3}$ of all 3-subsets of $X$, there are four orbits given by
\[\begin{array}{ll}
\mc{O}_1 = &  \{\{1, 2, 3\}, \{2, 3, 4\},  \{3, 4, 5\}, \{4, 5, 6\},  \{5, 6, 1\}, \{6, 1, 2\}\},\\
\mc{O}_2 = &  \{\{1, 2, 4\}, \{2, 3, 5\},  \{3, 4, 6\}, \{4, 5, 1\},  \{5, 6, 2\}, \{6, 1, 3\}\},\\
\mc{O}_3 = &  \{\{1, 2, 5\}, \{2, 3, 6\},  \{3, 4, 1\}, \{4, 5, 2\},  \{5, 6, 3\}, \{6, 1, 4\}\},\\
\mc{O}_4 = &  \{\{1, 3, 5\}, \{2, 4, 6\}\}.\end{array}\]
Suppose we take $P=X$ as point set and the union of two orbits $\mathcal{O}_1$ and $\mathcal{O}_4$ as block set. Then $(P, \mathcal{O}_1\cup\mathcal{O}_4)$ is  a PGD$(6,8,3, 4; 4,8)$ with circulant concurrence matrix
$NN^T=C[4, 2, 2, 0, 2, 2]$  and  Spec$(NN^T)=[12^1, 4^3, 0^2]$. 

\begin{remark}\label{Orbit}
\begin{enumerate}
\item The concurrence matrix indicates that \[\lambda_{xy} = \left \{\begin{array}{ll} 4 & \mbox{if } x=y,\\ 
0 & \mbox{if } \{x, y\}\in \{\{1, 4\}, \{2, 5\}, \{3, 6\}\},\\ 
2 & \mbox{otherwise.}\end{array}\right .\]
\item It is not a 2-design as there exists a pair $x,y\in X$ with $\lambda_{xy}=0$. 
\item It is a TD$_2(3, 2)$. 
\end{enumerate}
\end{remark}

\item[(b)] (\cite[Ex.8]{NO}) For the design in this example, we take the orbits $\mathcal{O}_1, \mathcal{O}_4$ and two copies of orbit $\mathcal{O}_2$ in the previous example with the same point set $X=\{1, 2, 3, 4, 5, 6\}$ and block size $k=3$. So the 3-sets in $\mathcal{O}_2$ are repeated blocks in this design $(X, \mathcal{O}_1\cup 2\cdot \mathcal{O}_2\cup \mathcal{O}_4)$. This is a PGD$(6, 20, 3, 10; 12, 18)$ with $\alpha=k(r-n)$ (so, $2$-$(6, 3, 4)$ design by \cite{N} (also cf. \cite[Remark 11]{NO})) and $NN^T = 10I+4(J-I)$, with $N=\left [N_1|N_2|N_2\right ]$ where $N_1$ being the incidence matrix $N$ of the previous example, and $N_2$ being the incidence matrix from the following:
\[N_2=\{\{1, 2, 4\}, \{3, 4, 6\}, \{2, 5, 6\}, \{2, 3, 5\}, \{1, 4, 5\}, \{1, 3, 6\}\}.\]

\item[(c)] (\cite[Ex.9]{NO}) Unlike the previous two examples, if we use group $G=\langle \sigma\rangle\cong \mathbb{Z}_5$ where $\sigma=(12345)\in S_6$, the cyclic permutation of length 5, then the 3-orbits of $X$ w.r.t. $G$ are:
\[\begin{array}{ll}
\mathcal{O}_1 = &  \{\{1, 2, 3\}, \{2, 3, 4\},  \{3, 4, 5\}, \{4, 5, 1\},  \{5, 1, 2\}\},\\
\mathcal{O}_2 = &  \{\{1, 2, 4\}, \{2, 3, 5\},  \{3, 4, 1\}, \{4, 5, 2\},  \{5, 1, 3\}\},\\
\mathcal{O}_3 = &  \{\{1, 2, 6\}, \{2, 3, 6\},  \{3, 4, 6\}, \{4, 5, 6\},  \{5, 1, 6\}\},\\
\mathcal{O}_4 = &  \{\{1, 3, 6\}, \{2, 4, 6\}, \{3, 5, 6\}, \{4, 1, 6\}, \{5, 2, 6\}\}.\end{array}\]
Since \[A_{3, 1} = \left [\begin{array}{cccccc}
3 & 3 & 3 & 3 & 3 & 0\\
3 & 3 & 3 & 3 & 3 & 0\\
2& 2 & 2 & 2 & 2 & 5\\
2& 2 & 2 & 2 & 2 & 5\end{array}\right ] \]
and $z=(1, 0, 0, 1)$ is a solution of $zA_{3,1}=5\mathbf{1}$,
by \cite[Theorem 6]{NO}, we have a PGD$(6, 10, 3, 5; 6, 9)$ $(X, \mathcal{O}_1\cup\mathcal{O}_4)$ with concurrence matrix 
$NN^T = 5 I + 2(J-I)$.
This is a $2$-$(6, 3, 2)$ design. 

\end{enumerate}

\subsection{An ad hoc construction on sets} 
 We conclude our survey of construction methods of PGDs with one elementary ad hoc construction given in \cite{OS}.
 \begin{example} \label{Ex3.1-OS}(Halved Ex.3.1 in \cite{OS}) Consider the tactical configuration $(P, \mathcal{B})$ defined by
\[P=G_1\cup G_2, \quad \mbox{where} \quad G_1=\{p_1, p_2, p_3\}, \ \ G_2= \{q_1, q_2, q_3\}; \]
\[\mathcal{B}=\{B_{ij}=S_i\cup T_j : S_i=G_1-\{p_i\}, T_j = G_2-\{q_j\}, i, j\in \{1, 2, 3\}\}\]
\[=\{B_{11}, B_{12}, B_{13}, B_{21}, B_{22}, B_{23}, B_{31}, B_{32},  B_{33}\}.\]

It is routine to verify that it is a PGD$(6, 9, 4, 6; 14, 17)$ with 
$NN^T=C[6, 4, 3, 4, 3, 4]$ and Spec$(NN^T)=[24^1, 3^4, 0^1]$. This PGD is isomorphic to a SPBIBD$(6, 9, 4, 6; 4, 3)$ based on a 2-class association scheme. The relation graphs of the underlying association scheme are $K_{3,3}$ and its complement.
 \end{example}

\section{Point-pair-concurrences in partial geometric designs} 
Neumaier \cite{N} determined the spectrum of the concurrence matrix of a PGD, and studied the classification problem of PGDs in terms of their parameters and spectra.  Van Dam and Spence \cite{vDS} used the spectra, when they produced a complete list of all PGDs with block size two and then went on to show that a PGD with block size three must be an SPBIBD, or a $2$-design. Recently, Lei, Qu and Shan \cite{LQS} took this result a step further and showed that a PGD with parameters $(v, b, 3, r; \alpha, \beta)$  is either: (i) a $2$-$(v, 3, \lambda)$ design, (ii) a TD$_{\lambda}(3, \frac{r}{\lambda})$, or (iii) each line of a generalized quadrangle that is repeated $\lambda$ times, where $\lambda = \frac{\beta}{2}+1$. Here a \emph{generalized quadrangle} is, by definition,  a partial geometry with $\tau = 1$. Lei, Qu and Shan also  considered PGDs with block size four for a particular concurrence type. We recall how the concurrences of a PGD have been used in the work of the above authors, and explore possibility of using concurrence matrices in the study of PGDs. 

\subsection{Point-pair-concurrence types} 

Given a block $B=\{x_1, x_2, \dots, x_k\}\in \mathcal{B}$ in a PGD $(P,\mathcal{B})$,  we define the \emph{concurrence type with respect to} $B$ by the $(k-1)$-tuple of integers
\[p_{Bx_1}:=[\lambda_{x_1x_2}, \lambda_{x_1x_3},\dots, \lambda_{x_1x_k}].\]
If PGD $\left(P,\mathcal{B}\right)$ has parameters $\left(v,b, k, r;\alpha, \beta\right)$, as we have seen in Lemma \ref{conc} and \ref{conc1}, we have  $$\sum\limits_{y\in B-\{x\}}\lambda_{xy}= \beta-r\quad \mbox{for every}\quad  x\in B\in\mathcal{B}.$$

While this sum is constant over all flags $(x, B)$, the $(k-1)$-tuple $p_{Bx}$ is not necessarily independent from the choice of $x$ and $B$ if $k>3$. The following example, provided by an anonymous reviewer, illustrates this.

\begin{example} \label{non-ciculant PGD} Consider the following PGD$(8, 8, 4, 4; 6, 10)$ with $P=\{1, 2, \dots, 8\}$ and $\mathcal{B}=\{B_1, B_2, \dots, B_8\}$ where $B_{i+1} = B_i^C$, for $i=1, 3, 5, 7$,  and 
\[\begin{array}{|c|ccccccc|}
  i & 1 & & 3 & & 5 &  & 7  \\ \hline
B_i& \{1, 2, 3, 4\} & & \{1, 2, 5, 6\} & & \{1, 3, 5, 7\} & & \{1, 4, 5, 8\} \end{array}\ .\]
 Its concurrence matrix may be described as  
 \[NN^T = \left [\begin{array}{c|c} 2I+2J & 2I + J\\ \hline 2I+J & 2I + 2J \end{array}\right ],\] where $I$ and $J$ are the $4\times 4$ identity matrix and all-ones matrix, respectively. In this PGD, the concurrence type  of point 1 with respect to block $B_1$ and with respect to $B_3$ are different: that is, $p_{B_11}=[\lambda_{12}, \lambda_{13}, \lambda_{14}]=[2, 2, 2]$ and $p_{B_31} = [\lambda_{15}, \lambda_{12}, \lambda_{16}]=[3, 2, 1]$.
 
We observe that this $NN^T$ is not circulant and Spec$(NN^T)=[16^1, 4^4, 0^3]$.\end{example}

\subsection{PGDs with block sizes 2 and 3} 
The classification from van Dam and Spence \cite{vDS} concerning PGDs with block size two is stated: 
\begin{theorem}\cite{vDS}
Let $D = (P, \mathcal{B})$ be a PGD with parameters $(v, b, 2, r; \alpha, \beta)$.  Then, $D$ is of the form $D_1 \otimes J_{1, \beta+1}$ where $D_1$ is either a TD$_1(2,\frac{r}{\beta+1})$, or the design consisting of all $\binom{v}{2}$ pairs of points from $P$.
\end{theorem}

\begin{proof}
See Proposition 6 in \cite{vDS}.
\end{proof}

 The results from \cite{LQS} concerning PGDs with block size three are formally stated as follows:
\begin{lemma}\cite{LQS}{\label{concprofsizethree}}
Let $D = (P, \mathcal{B})$ be a PGD$(v, b, 3, r; \alpha, \beta)$.  Then, for any block $B = \{x_1, x_2, x_3\} \in \mathcal{B}$, we have:
$$\lambda_{x_1x_2} = \lambda_{x_1x_3} = \lambda_{x_2x_3}.$$
That is, the concurrence type $p_{Bx_1}$ is $[\lambda, \lambda, \lambda]$ for some positive integer $\lambda$.
\label{inprofblockthree}
\end{lemma}

\begin{proof}
Consider the flags $(x_1,B), (x_2,B)$, and $(x_3,B)$. Then we have:
$$\beta = \lambda_{x_1x_1} + \lambda_{x_1x_2} + \lambda_{x_1x_3}
 = \lambda_{x_2x_1} + \lambda_{x_2x_2} + \lambda_{x_2x_3}
 = \lambda_{x_3x_1} + \lambda_{x_3x_2} + \lambda_{x_3x_3}.$$
Using the fact that $\lambda_{x_ix_i} = r$ and $\lambda_{x_ix_j} = \lambda_{x_jx_i}$ for any points $x_i$ and $x_j$, we get that $\lambda_{x_1x_2} = \lambda_{x_1x_3} = \lambda_{x_2x_3}$ by solving this system of linear equations.
\end{proof}

\begin{theorem}\cite{LQS}
Let $D = (P, \mathcal{B})$ be a PGD$(v, b, 3, r; \alpha, \beta)$ and $\lambda = \frac{\beta}{2}+1$.  Then one of the following holds:
\begin{enumerate}
    \item D is a $2$-$(v, 3, \lambda)$ design, 
    \item D is a TD$_ {\frac{r}{\lambda}}(3, \lambda)$, or
    \item D has the form $D_1 \otimes J_{1, \lambda}$ where $D_1$ is a generalized quadrangle of order $(2,2)$ or $(2,4)$.
\end{enumerate}
\end{theorem}
\begin{proof}
We refer the reader to \cite{LQS}.
\end{proof}

We note that a PGD in Part 2 or 3 is isomorphic to an SPBIBD (cf. \cite[Prop.7]{vDS}).
\subsection{$[\lambda, 1, 1]$-type PGDs with block size 4}
For the PGDs with block size 4, the following was also shown in \cite{LQS}:
\begin{lemma}\cite{LQS}
Let $D = (P, \mathcal{B})$ be a PGD with parameters $(v, b, 4, r; \alpha, \beta)$.  Then, for any block $B = \{x_1, x_2, x_3, x_4\} \in \mathcal{B}$, 
$$\lambda_{x_1x_2} = \lambda_{x_3x_4}, \hspace{1cm} \lambda_{x_1x_3} = \lambda_{x_2x_4}, \hspace{1cm} \lambda_{x_1x_4} = \lambda_{x_2x_3}.$$
\label{inprofblockfour}
\end{lemma}

\begin{proof}
The proof is similar to the proof of Lemma~\ref{concprofsizethree}.
\end{proof}

The following results for a PGD with block size 4 is reported in \cite{LQS}.
\begin{theorem}\cite{LQS}
Let $D = (P, \mathcal{B})$ be a PGD$(v, b, 4, r; \alpha, \beta)$. Suppose $D$ has a property that $p_{Bx} =[\lambda, 1, 1]$ for any block $B\in \mathcal{B}$.  Then, exactly one of the following holds:
\begin{enumerate}
    \item $\lambda = 1$ and $D$ is a partial geometry.
    \item $\lambda > 1$ and $\alpha \in \{2, 4\}$.
    \begin{enumerate}
        \item Suppose $\alpha=4$. 
        \begin{enumerate}
            \item $r=\lambda$ and $D$ is of the form $2$-$(r+1, 2, 1) \otimes J_{2,1}$.
            \item $r > \lambda = 2$ and $D$ gives rise to a $3$-class association scheme with intersection matrices, $B_1, B_2$ and $B_3$, respectively, given by
            \[\left [ \begin{array}{cccc} 0 & 1 & 0 & 0\\ \frac{r}{2} & 0 & 0 & \frac{r}{2}\\ 0 & 0 & \frac{r}{2} & 0\\
            0 & \frac{r-2}{2} & 0 & 0\end{array}\right ],\ \
            \left [ \begin{array}{cccc} 0 & 0 & 1 & 0\\ 0 & 0 &  \frac{r}{2} & 0 \\ 2r & 2r & r & 2r\\  
            0 & 0 & \frac{r-2}{2} & 0\end{array}\right ],\ \
            \left [ \begin{array}{cccc} 0 & 0 & 0 & 1\\ 0 & \frac{r-2}{2} & 0 & 0 \\ 0 & 0 & \frac{r-2}{2} &  0\\  \frac{r-2}{2} & 0 
            & 0 & \frac{r-4}{2}\end{array}\right ].\]
        \end{enumerate}
        \item Suppose $\alpha = 2$.
        \begin{enumerate}
            \item $r=\lambda$ and $D$ is of the form TD$_1(2, r) \otimes J_{2,1}$.
            \item $r>\lambda = 2$ and $D$ has parameters $(32, 48, 4, 6; 2, 10)$. A specific construction of a PGD$(32, 48, 4, 6; 2, 10)$ is given in Section 4 of \cite{LQS}.

        \end{enumerate}
    \end{enumerate}
\end{enumerate}
\end{theorem}

\begin{proof}
We refer the reader to \cite{LQS}.
\end{proof}

\subsection{PGDs of block size 4}
We consider three examples of PGDs with block size $4$ that have other concurrence types with respect to a block. These examples will indicate that the complete classification of PGDs with block 4 requires a lot more work.
\begin{example}{\label{28-2.10}}
In \cite[Example 3.14]{O}, Olmez showed that a PGD $D=(P,\mathcal{B})$, defined by $P=\mathbb{Z}_2 \times \mathbb{Z}_4$ and $\mathcal{B} = \{g+S:g \in \mathbb{Z}_2 \times \mathbb{Z}_4\}$ exists where $S = \{(0,0),(0,2),(0,3),(1,1)\}$ is a PGDS in $\mathbb{Z}_2 \times \mathbb{Z}_4$. $D$ is a PGD$(8,8,4,4;6,10)$ and is isomorphic to the PGD discussed in Example \ref{PGDS1} by the mapping from $\{(0,0),(0,1),(0,2),(0,3),(1,0),(1,1),(1,2),(1,3)\}$ to $\{1, i, j, k, -1, -i, -j, -k\}$, in order.  The   concurrence matrix $NN^T$ may be described as:
\[NN^T = C[4, 2, 2, 2, 0, 2, 2, 2] \ \mbox{ and }\ 
\lambda_{xy} = \left \{ \begin{array}{ll}
      0 & \text{if } \, x+(1,0)=y \,\,\text{in }\,\mathbb{Z}_2 \times \mathbb{Z}_4,\\ 
      2 & \text{neither }x=y,\mbox{ nor } x+(1,0)= y.\\ 
   \end{array}\right .\]
Thus, in this PGD,  $p_{Bx} =[2,2,2]$ with respect to any block $B\in\mathcal{B}$. 
\end{example}
Note that all $2$-$(v, 4, 2)$ designs also have concurrence type $[2,2,2]$ with respective to any block, but the PGDs in the previous example are not $2$-designs.  This is because not all pairs of points occur in a block together. Thus, we have a PGD with concurrence type $[2,2,2]$  with respect to any block that is not a $2$-design nor a $\mathcal{D}\otimes J_{1, 2}$ with $\mathcal{D}$ being a PGD of type $[1, 1, 1]$. 

\begin{example}{\cite{O}}
Let $P = \mathbb{Z}_8$ and $S = \{0, 1, 4, 5\}$.  Then $(\mathbb{Z}_8,\mathcal{B})$ with
$$\mathcal{B} = \{ g + S: g \in \mathbb{Z}_8\}$$
is a PGD$(8, 8, 4, 4; 4, 12)$. The concurrences are: 
\[NN^T=C[4, 2, 0, 2, 4, 2, 0, 2] \quad \mbox{ and } \quad 
\lambda_{xy} = \begin{cases} 
      4 & \text{if} \,\,\, x=y \,\,\, or \,\,\, \left| x-y \right| \equiv 4\pmod{8}\\
      0 & \text{if} \,\,\, \left| x-y \right| \equiv 2\pmod{8} \\
      2 & \text{otherwise.}
   \end{cases}\]
The concurrence type with respect to any block of this PGD is $[4,2,2]$. This PGD is not isomorphic to a PGD $\mathcal{D}\otimes J_{1, 2}$ with $\mathcal{D}$ being a PGD of type $[2, 1, 1]$. 
\end{example}

\begin{example} (Example \ref{Ex3.1-OS}--revisited.)  
In Example \ref{Ex3.1-OS}, we have seen a PGD$(6,9,4,6; 14,17)$ with
\[NN^T=C[6, 4, 3, 4, 3, 4]\quad \mbox{ and }\quad 
\lambda_{xy} = \begin{cases} 
      3 & \text{if either} \,\,\,\, x, y \in G_1 \,\,\,\, \text{or} \,\,\,\, x, y \in G_2 \\
      4 & \text{otherwise.}
   \end{cases}\]
The concurrence type with respect to any block of this PGD is  $[4,4,3]$.
\end{example}

\subsection{PGDs having circulant concurrence matrices}
We notice that the concurrence matrices of many PGDs are circulant. We close this section by introducing an infinite family of PGDs whose concurrence matrices are circulant.

Consider the $2$-$(v, 2, 1)$ design $D=(P, \mathcal{B})$ where $P=\{1, 2, \dots, v\}$ and $\mathcal{B}$ is the set of all two-element subsets of $P$, with $b=v(v-1)/2$ and $r=v-1$. The incidence matrix $N_{(v)}$ for this design can be described as the following. 
\[N_{(3)}=\left [\begin{array}{ccc} 1 & 1 & 0\\ 1 & 0 & 1\\ 0 & 1 & 1\end{array}\right ]\quad \mbox{ and }\quad 
N_{(v)} = \left [\begin{array}{c|c}
1\ 1\ \cdots \ 1 & 0\ 0\ \cdots \ 0\\ \hline
& \\
I_{v-1} & N_{(v-1)} \\
& \\
\end{array}\right ]\quad \mbox{for } v\ge 4.\]
Then $N_{(v)}{N_{(v)}}^T=C[v-1, \underbrace{1, 1, \dots, 1}_{v-1}]$. Moreover, for $N=N_{(v)}\otimes J_{2,1}$, which is the incidence matrix of $D\otimes J_{2,1}$ with parameters $(2v, {v\choose 2}, 4, v-1; 4, 2v)$ satisfies
\[NN^T=[N_{(v)} \otimes J_{2,1}][N_{(v)} \otimes J_{2,1}]^T = C[v-1, \underbrace{1, 1, \dots, 1}_{v-1}, v-1, \underbrace{1, 1, \dots, 1}_{v-1}].\]
Hence we have:

\begin{prop}
For any integer $v \geq 3$, there is a PGD$(2v, {v\choose 2}, 4, v-1; 4, 2v)$ whose concurrence matrix is circulant.
The spectrum of this PGD is $[4(v-1)^1, (2v-4)^{v-1}, 0^v]$. 
\end{prop}
\begin{proof}
It follows from $NN^T=N_{(v)}{N_{(v)}}^T\otimes J_{2, 1}$, and the rest of the proof is straightforward.
\end{proof}

In general, we have the following.
\begin{prop}
Suppose $N$ be the incidence matrix of a PGD $\mathcal{D}$ whose concurrence matrix is circulant. Then for any $m,n\in\mathbb{N}$, the concurrence matrix of $\mathcal{D}\otimes J_{m,n}$ is circulant.
\end{prop} 
\begin{proof}
It is straightforward from the fact that $[C\otimes J_{m,n}][C\otimes J_{m,n}]^T$ is circulant for any circulant matrix $C$. 
\end{proof}

As a consequence,  for instance, we have the following.
\begin{cor}
Suppose $\mathcal{D}$ is a 2-design. Then the concurrence matrix of the PGD
 $\mathcal{D}\otimes J_{1,n}$ is circulant for any $n\in\mathbb{N}$.
\end{cor}
\begin{proof}
It is immediate from the previous proposition.
\end{proof}

%

\section{Classification of PGDs having circulant concurrence matrices}
We have observed that the concurrence matrices of many PGDs are circulant. Motivated by this observation, in this section, we explore small PGDs which have circulant concurrence matrices.

It is well known that the eigenvalues of a real $v \times v$
circulant matrix $C=C[c_0, c_1, \dots, c_{v-1}]$
are given by $f(\omega^k)$ for $k = 0, 1, \hdots, v-1$ where $\omega$ is a primitive $v$th-root of unity  
and 
\[f(\lambda) = c_0 + c_1 \lambda + c_2 \lambda ^2 + \cdots +c_{v-2} \lambda ^{v-2} + c_{v-1} \lambda ^{v-1}.\]  
Note that $C$ is symmetric if and only if $c_{v-i}=c_i$ for $i = 1, 2, \dots, \lfloor{\frac{v}{2}}\rfloor$.
We will denote a symmetric $v\times v$ circulant matrix $C$ by:
$$\begin{cases} 
      C[c_0, c_1, \dots, c_{\frac{v}{2}}, c_{\frac{v}{2}-1}, \dots, c_1] & \text{if $v$ is even} \\
       C[c_0, c_1, \dots, c_{\frac{v-1}{2}}, c_{\frac{v-1}{2}}, c_{\frac{v-1}{2}-1}, \dots, c_1] & \text{if $v$ is odd}
   \end{cases}$$
 or simply by 
 \[C[c_0, c_1, \dots, c_{\lfloor\frac{v}{2}\rfloor}].\]

Our aim is to find all circulant matrices of order up to 12 that may be realized as feasible concurrence matrices of PGDs. We then determine whether each of these matrices is actually realized as the concurrence matrix of a PGD by describing corresponding incidence matrix (getting computer aid if necessary). If it is realized, then we describe and classify the PGD: that is, determine if it is isomorphic to any known PGDs.
We do this for each fixed order $v$, for $v=6, 8, 9, 10$ and $12$, and our study is based on the following facts. 

\begin{remark}\label{Spec1}
Let $C$ be a $v\times v$ symmetric circulant matrix $C[c_0, c_1, \dots, c_{\lfloor\frac{v}{2}\rfloor}].$
\begin{enumerate}
\item The eigenvalues $\theta_j=f(\omega^j), j=0, 1, 2, \dots, v-1$ for $C$ are given by
\[\theta_j = \left\{ \begin{array}{ll}
c_0 +2 c_1 \mathcal{R}\omega^{1j} + 2c_2 \mathcal{R}\omega^{2j} + \cdots +2c_{v/2-1} \mathcal{R}\omega^{(v/2-1)j} + c_{v/2} \mathcal{R}\omega^{(v/2)j} & \mbox{for } v \mbox{ even}\\
c_0 +2 c_1 \mathcal{R}\omega^{1j} + 2c_2 \mathcal{R}\omega^{2j} + \cdots +2c_{(v-1)/2} \mathcal{R}\omega^{((v-1)/2)j} & \mbox{for } v \mbox{ odd,}\end{array}\right .\]  where $\omega$ is a primitive $v$th-root of unity 
(so, $\omega^v-1 = 0$ and $\omega^k-1 \neq 0$ for every positive integer $k < n$) and 
$\mathcal{R}\omega^{kj} = \cos(2kj\pi/v)$  (cf. \cite{Da}).

\item Suppose $C$ is a concurrence matrix of a PGD$(v, b, k, r; \alpha, \beta, n)$.
 Then, in Lemma \ref{Spec},  \[ c_0= r, \quad Spec(C) =[kr^1, n^{\sigma}, 0^{v-1-\sigma}]\quad \mbox{where} \quad \sigma=r(v-k)/n.\]
\end{enumerate}
\end{remark}
We note that we can reduce some parameters of a PGD$(v, b, k, r; \alpha, \beta, n)$ according to the following relations derived in Proposition \ref{parameters}: \[b=vr/k, \ \ \alpha=k(kr-n)/v, \ \ \beta=n+\alpha = n+k(kr-n)/v.\] Utilizing these relations, we will sometimes denote PGD$(v, b, k, r; \alpha, \beta, n)$ by PGD$(v, k, r; n)$ with another essential parameter $\sigma= r(v-k)/n$.
Working under the assumption that the circulant matrix is feasible as a concurrence matrix of a PGD, it is necessarily $c_0 =r$.  
The rest of the entries of the putative circulant matrix are determined by solving a system of equations obtained from each given set of possible values for $k, n, \sigma$ and $r$ for given $v$. 

In the rest of this section, we summarize our findings according to the order of the matrices from $v=6$ to $v=12$, except for primes $v=7$ and $11$. (We note that 
for a prime $p$, the eigenvalue $f(\omega^j)$ for any $j\in \{0, 1, 2,\dots, p-1\}$  ($\omega$ being the $p$th roots of unity) is integer only if $c_1=c_2=\cdots =c_{p-1}$. Therefore, $C[r, n, n, \dots]$ appears as the only feasible circulant matrix which is associated with a 2-design if there is any nontrivial PGD.)\footnote{We recall that if $\omega^v=1$ and $d|v$, say $v=dl$, then $(\omega^d-1)|(\omega^v-1)$ and $\omega^d + \omega^{2d}+\cdots + \omega^{(l-1)d}\in \mathbb{Z}$, and so, $f(\omega^d)\in\mathbb{Z}$.} 

\subsection{PGDs of order 6} 

We classify all $6\times 6$ circulant matrices $C[c_0, c_1, c_2, c_3, c_2, c_1]$ with nonnegative 
integers $c_0, c_1, c_2, c_3$, that are realized as the concurrence matrices of PGDs, and then describe these PGDs.
According to Remark \ref{Spec1}, for a $6\times 6$ symmetric circulant matrix $C[c_0, c_1, c_2, c_3]$,  the eigenvalues of $C$ must be: 
$$\theta_j =  c_0 + 2c_1 \cos\left(\frac{j \pi}{3}\right) + 2c_2 \cos\left(\frac{2j \pi}{3}\right) + c_3 \cos\left(j \pi\right).$$
Since the values for $\cos(\frac{j \pi}{3})$, $\cos(\frac{2j \pi}{3})$, and $\cos(j \pi)$ for each $j \in \{0,1,2,3,4,5\}$ lie in $\{ \pm \frac{1}{2}, \pm 1\}$,  $C$ is to be $NN^T$ for a PGD$(v, k,r;n)$ with incidence matrix $N$, all feasible values $c_1, c_2, c_3$ may be determined by considering all appropriate combinations of the spectra coming from the following table:

\begin{center}
\begin{tabular}{ |c|c|c|c| } 
 \hline
 $j$ & $\theta_j$ & Expected value(s) for $\theta_j$ & Multiplicity \\
 \hline
 0 & $c_0 + 2c_1 + 2c_2 + c_3$ & $kr$ & 1 \\ 
 $1$, $5$ & $c_0+c_1-c_2-c_3$ & either $0$ or $n$ & 2 \\ 
 $2$, $4$ & $c_0-c_1-c_2+c_3$ & either $0$ or $n$ & 2 \\ 
 $3$ & $c_0-2c_1+2c_2-c_3$ & either $0$ or $n$ & 1 \\ 
 \hline
\end{tabular}
\end{center}

Then, a PGD$\left (6, \frac{6r}{k}, k, r; \frac{k(kr-n)}{6}, n+\frac{k(kr-n)}{6}\right )$ (using the equations from Prop. \ref{parameters}) will have spectrum 
$$[kr^1, n^{\sigma}, 0^{5-\sigma}] \,\,\,\,\,\, \text{where} \,\,\,\,\, \sigma = \frac{r(6-k)}{n}.$$

The feasible solutions $(c_0, c_1, c_2, c_3)$ are listed below according to the choices of multiplicity $\sigma$ of eigenvalue $n$ for the case $k=3$ and $k=4$.

\flushleft 1. \textbf{Case when $\underline{k=3}$}: For $k=3$, we have that $2 \leq  \sigma = \frac{3r}{n} \leq 4.$ (Note that $\sigma = 1$ is not a possibility as it will imply that $n=kr$ and $\alpha = 0.$)
\begin{enumerate}
        \item[(a)] If $\sigma = 2 \,\,(n = \frac{3r}{2})$, then $\theta_3 = 0$ and either (i) $\theta_{1}=\theta_{5} = n$ and $\theta_{2}=\theta_{4} = 0$ or (ii) $\theta_{1}=\theta_{5} = 0$ and $\theta_{2}=\theta_{4} = n.$  By solving the system for the $c_i$ values, we get:
        $$(c_0,c_1,c_2,c_3)= \begin{cases} 
      (r,\frac{3r}{4},\frac{r}{4},0) & \,\,\,\,\,\,\,\,\,\ \text{for case (i)} \\
      (r,\frac{r}{4},\frac{r}{4},r) & \,\,\,\,\,\,\,\,\,\ \text{for case (ii)}
   \end{cases}$$
    
    For $l \in \mathbb{N}$, 
         $ C= C[4l,3l,l,0,l,3l] \mbox{ and }  C[4l,l,l,4l,l,l]$ appear to be possible concurrence matrices for a PGD with $k=3, \sigma=2, r=4l, n=6l$. However, it is easy to verify that none of these matrices are realized as a concurrence matrix of a PGD with the parameters.
         
    \item[(b)] If $\sigma = 3 \,\,(n = r)$, then $\theta_3 = n$ and either (i) $\theta_{1}=\theta_{5} = n$ and $\theta_{2}=\theta_{4} = 0$ or (ii) $\theta_{1}=\theta_{5} = 0$ and $\theta_{2}=\theta_{4} = n.$  By solving the system for the $c_i$ values, we get:
        $$(c_0,c_1,c_2,c_3)= \begin{cases} 
      (r,\frac{r}{2},\frac{r}{2},0) & \,\,\,\,\,\,\,\,\,\ \text{for case (i)} \\
      (r,\frac{r}{6},\frac{r}{2},\frac{2r}{3}) & \,\,\,\,\,\,\,\,\,\ \text{for case (ii)}
   \end{cases}$$

    For $l \in \mathbb{N}$,  
         $\mbox{(i)}\ C[2l,l,l,0,l,l]\ \mbox{and } \mbox{(ii)}\ C[6l,l,3l,4l,3l,l]$ appear to be possible concurrence matrices for a PGD with $k=3, \sigma=3, n=r$. For each $l\in\mathbb{N}$, the circulant matrix $C[2l,l,l,0,l,l]$ in (i) is realized as the concurrence matrix of a PGD$(6, 4l, 3, 2l; 2l, 4l)$. We demonstrate these PGDs below. However, it is easy to verify that none of the matrices in (ii) can be realized as $NN^T$ for a PGD. 

\begin{example} $C = C[2, 1, 1, 0, 1, 1]$ is realized as $NN^T$  for the PGD$(6, 4, 3, 2; 2, 4)$\  $(P, \mathcal{B})$  where  $P=\{1, 2, 3, 4, 5, 6\}$ and 
$ \mathcal{B} = \{\{1,2,3\}, \{1,5,6\}, \{2,4,6\}, \{3,4,5\}\}.$
Note that this PGD is isomorphic to a TD$_1(3,2)$.
\end{example}

\begin{example} $C = C[4, 2, 2, 0, 2, 2]$ is realized as $NN^T$ for the PGD$(6, 8, 3, 4; 4, 8)$ whose block set is
\[\mathcal{B} = \left \{\begin{array}{l} \{1,2,3\}, \{1,2,6\}, \{1,3,5\}, \{1,5,6\},\\
 \{2,3,4\}, \{2,4,6\}, \{3,4,5\}, \{4,5,6\}\end{array}\right \}.\]
Note that this PGD is isomorphic to a TD$_2(3,2)$.
\end{example}

\begin{remark} 
\begin{enumerate}
\item[$(1)$] It is easy to see that each PGD$(6,4l,3,2l; 2l, 4l)$ with $l\in\mathbb{N}$ is isomorphic to a TD$_l(3,2)$.   (Also see Remark \ref{Orbit}.)

\item[$(2)$] These PGDs with block size less than 30 are all accounted for in the list of small PGDs in \cite{vDS}; namely, for $l=2, 3, 4, 5, 6, 7$ they are listed with labels N1, N3, N13, N24, N43, and N76, respectively. \end{enumerate}
\end{remark}

\item[(c)] If $\sigma = 4 \,\,(n = \frac{3r}{4})$, then $\theta_3 = 0$ and $\theta_{1}=\theta_{2}=\theta_{4}=\theta_{5} = n$. The only solution here is  $(c_0, c_1, c_2, c_3) = (4l, 2l, l, 2l)$ for $l \in \mathbb{N}$.  However, there does not exists a PGD under required parameters with $NN^T = C[4l,2l,l,2l,l,2l]$ for any $l$.
    \end{enumerate}

2. \textbf{Case when $\underline{k=4}$}: For $k=4$, we have that $1 \leq  \sigma = \frac{2r}{n} \leq 4.$\begin{enumerate}
    \item[(a)] If $\sigma=1 \ \ (n=2r)$, then $\theta_3 = n$ and  $\theta_{1}=\theta_{2} = \theta_{4}=\theta_{5} = 0$. The solution
$(c_0, c_1, c_2, c_3)=(r, \frac{r}{3}, r, \frac{r}{3})$ to the corresponding system of equations from Table 1 yields symmetric circulant matrices $C=C[6l, 2l, 6l, 2l, 6l, 2l]$ for $l\ge 1$. For each $l$, $C$ could be a concurrence matrix for a putative PGD$(6, 9l, 4, 6l; 8l, 20l)$. However, it is shown that there does not exist a PGD$(6, 9l, 4, 6l; 8l, 20l)$ for any $l\ge 1$.

    \item[(b)] If $\sigma = 2 \,\,(n = r)$, then $\theta_3 = 0$ and either (i) $\theta_{1}=\theta_{5} = n$ and $\theta_{2}=\theta_{4} = 0$ or (ii) $\theta_{1}=\theta_{5} = 0$ and $\theta_{2}=\theta_{4} = n.$ By solving the system for the $c_i$ values, we get:
        $$(c_0,c_1,c_2,c_3)= \begin{cases} 
      (r,\frac{5r}{6},\frac{r}{2},\frac{r}{3}) & \,\,\,\,\,\,\,\,\,\ \text{for case (i)} \\
      (r,\frac{r}{2},\frac{r}{2},r) & \,\,\,\,\,\,\,\,\,\ \text{for case (ii)}
   \end{cases}$$
Thus, for $l\in\mathbb{N}$,\ $\mbox{(i)}\ C= C[6l,5l,3l,2l,3l,5l]\ \ \mbox{and}\ \ \mbox{(ii)}\  C[2l,l,l,2l,l,l]$ appear to be possible concurrence matrices for PGDs. 

We can verify that $C$ in the case of (i) is not suitable for $NN^T$ for any PGD. 
However, for case (ii), $C[2l,l,l,2l,l,l]$ for each $l\in \mathbb{N}$ is realized as $NN^T$ for a PGD($6,3l,4,2l;4l,6l)$ (with $\sigma=2, r=n=2l$). Furthermore, for each $l$, PGD($6,3l,4,2l;4l,6l)$ is unique and is isomorphic to a $2$-$(3, 2, 1) \times J_{2,l}$. In particular, 
 $C = C[2, 1, 1, 2, 1, 1]$ is realized as $NN^T$
for the  PGD$(6, 3, 4, 2; 4, 6)$ with block set 
$\mathcal{B} = \{\{1,2,4,5\}, \{1,3,4,6\}, \{2,3,5,6\}\}.$

\item[(c)] If $\sigma = 3 \,\,(n = \frac{2r}{3})$, then $\theta_3 = n$ and either (i) $\theta_{1}=\theta_{5} = n$ and $\theta_{2}=\theta_{4} = 0$ or (ii) $\theta_{1}=\theta_{5} = 0$ and $\theta_{2}=\theta_{4} = n.$  By solving the system (i) or (ii), we obtain the same solution set for the $c_i$ values, namely, $(c_0, c_1, c_2, c_3) = (r, \tfrac{2r}{3}, \tfrac{2r}{3}, \tfrac{r}{3})$, and thus,  we obtain the set of putative circulant matrices, $C[6l,4l,4l,2l,4l,4l].$ However, none of the circulant matrices can be realized as a concurrence matrix.

\item[(d)] If $\sigma = 4 \,\,(n = \frac{r}{2})$, then $\theta_3 = 0$ and $\theta_{1}=\theta_{2}=\theta_{4}=\theta_{5} = n$. The only solution here is $(c_0, c_1, c_2, c_3) = (6l, 4l,3l, 4l)$ for $l \in \mathbb{N}$ (with $\sigma=4, r=6l$ and $n=3l$).  For each $l\in\mathbb{N}$, $C[6l, 4l,3l, 4l, 3l, 4l]$ is realized as $NN^T$ of a PGD($6, 9l, 4, 6l; 14l, 17l$).    

As $\sigma=v-2$, all the PGDs in this family belong to the family of SPBIBDs. Furthermore, it can be shown with a simple computer algorithm that each of these PGDs is unique.
In the case of $l=1$, $C =C[6, 4, 3, 4, 3, 4]$ is the concurrence matrix of a PGD($6,9,4,6;14,17$) which is isomorphic to the PGD described in Example \ref{Ex3.1-OS}. 

\end{enumerate}
Note that in the above classification, since PGDs with block size 4 has been studied, those with block size 2 as their complementary conjugates have also been studied. They are improper PGDs under the condition $3\le k\le v-3$; for this reason, none of these PGDs are discussed in \cite{vDS}.
We conclude our analysis and classification of PGDs of order 6 having circulant concurrence matrices with the following proposition. 

\begin{prop}
There is an infinite family of PGDs of order 6 with block size 3
and two infinite classes of PGDs with block size 4 whose concurrence matrices are circulant as in the table
\label{O6}
\[\begin{array}{|c|c|cc|cc|} \hline
\mbox{Case} \#   &   C= NN^T             &   n^{\sigma} &\mbox{PGD}(v, b, k, r; \alpha, \beta)   & {Remark} & \\ \hline
1-(b)  & C[2l, l, l, 0, l,  l]        &   (2l)^3 &  (6, 4l, 3, 2l; 2l, 4l) &   TD_l(3,2)  & \\
2-(a) &  C[2l, l, l, 2l, l, l]        &    (2l)^2 &  (6, 3l, 4, 2l; 4l, 6l)  & M\otimes J_{2, l}&\\
2-(c) & C[6l, 4l, 3l, 4l, 3l, 4l] &     (3l)^4  & (6, 9l, 4, 6l; 14l, 17l) & N\otimes J_{1, l} & \\ \hline
\end{array}\] where $M$ is the incidence matrix of a $2$-$(3, 2, 1)$ design, and $N$ is the incidence matrix of a PGD$(6, 9, 4, 6; 14, 17)$ given in Example \ref{Ex3.1-OS}.
\end{prop}

\begin{proof} We have seen that the circulant matrix $C[2l, l, l, 0, l, l]$ for the case 1-(b) is realized as the concurrence matrix of a PGD$(6, 4l, 3, 2l; 2l, 4l)$ in the above analysis. Therefore, it is sufficient to assert that  the concurrence matrix of a PGD$(6, 4l, 3, 2l; 2l, 4l)$  is $C[2l, l, l, 0, l, l]$ to complete the proof. We omit the details.
\end{proof}


\subsection{PGDs of order 8} 
When $v=8$, $C=C[c_0,c_1,c_2,c_3,c_4, c_3,c_2,c_1]$ and has eigenvalues
\[\theta_j=c_0+2c_1\cos(j\pi /4)+2c_2 \cos(j\pi/2)+2c_3\cos(3j\pi/4) +c_4\cos(j\pi),\]
where the values $\cos(j\cdot \pi /4),  \cos(j\cdot \pi /2), \cos(3j\pi/4)$, and $\cos(j\cdot \pi)$ for each $j=1, 2, \dots, 7$ lie in $\{\pm \frac{\sqrt{2}}{2}, 0, \pm 1\}$. As we have done for the case with $v=6$, we derive sets of system of equations under the given constraints for the values of $c_1, c_2, c_3$ and $c_4$ and the spectrum of $C$ for every feasible combination of $k$ and $\sigma$ referring to the following table.
\[\begin{array}{|r|c|c|c|c|}\hline
j & \theta_j & \mbox{Expected values for }\theta_j & \mbox{multiplicity}& \mbox{remark}\\ \hline
0 & c_0+2c_1+2c_2+2c_3+c_4 & kr & 1 & \\
1,3, 5,7 & c_0-c_4 & \mbox{either } 0 \mbox{ or } n & 4 &c_1=c_3\\
2, 6 & c_0-2c_2+c_4 & \mbox{either } 0 \mbox{ or } n & 2 & \\
4 & c_0-2c_1+2c_2-2c_3+c_4 & \mbox{either } 0 \mbox{ or } n & 1&\\ \hline \end{array}\]

We note that $c_1=c_3$ can be seen by looking at the case when $j=1$. Namely, 
\[f(\omega) = c_0 + c_1\cos(\tfrac{\pi}{4}) + c_2\cos(\tfrac{\pi}{2}) + c_3\cos(\tfrac{3\pi}{4}) + c_4\cos(\pi) + c_3\cos(\tfrac{5\pi}{4}) + c_2\cos(\tfrac{3\pi}{2}) +c_1\cos(\tfrac{7\pi}{4})\]
\[\qquad\qquad = \quad c_0 +\sqrt{2}c_1-\sqrt{2}c_3-c_4.\]
Since the eigenvalues must be integral, it must be the case that $c_1-c_3=0$. Similarly, $c_1-c_3=0$ for the cases $j=3, 5, 7$.

Suppose the $8\times 8$ symmetric circulant matrix $C$ has Spec$(C)=[kr^1, n^{\sigma}, 0^{7-\sigma}]$, where $n=r(8-k)/\sigma$ for $\sigma\in \{1, 2, 3, 4, 5, 6\}$. Then the suitable entries $c_1, c_2, c_3, c_4$ for putative $C$ must appear in a row of the following table.

\[\begin{array}{|c|c||c|c|c||c|c|c|c|}\hline
 \sigma & n & \theta_1, \theta_3, \theta_5, \theta_7 & \theta_2, \theta_6 & \theta_4 & c_1 & c_2 & c_3 & c_4    \\ \hline\hline
1 & r(8-k) & 0 & 0 & n & r-\tfrac{n}{4} & r & r-\tfrac{n}{4} & r  \\ \hline
  2 &\frac12r(8-k)& 0 & n & 0 & r-\frac{n}{4} & r-\frac{n}{2} & r-\frac{n}{4} & r   \\  \hline
 3 & \frac13r(8-k)& 0 & n & n &  r-\frac{n}{2} & r-\frac{n}{2} & r-\frac{n}{2} & r   \\ \hline
 4 & \frac14r(8-k)& n & 0 & 0 &  r-\frac{n}{2} & r-\frac{n}{2} & r-\frac{n}{2} & r-n   \\ \hline
 5 & \frac15r(8-k)& n & 0 & n &  r-\frac{3n}{4} & r-\frac{n}{2} &r-\frac{3n}{4} & r-n   \\ \hline
 6 & \frac16r(8-k)& n & n & 0 &  r-\frac{3n}{4} & r-n & r-\frac{3n}{4}  & r -n   \\ \hline
\end{array}\]

Furthermore, for each $k=3$ and $k=4$, we have the following putative matrices as concurrence matrices of PGDs.

\begin{enumerate}
\item For $k=3$, each of the cases for $\sigma=1, 2$ and $4$ cannot yield a feasible value of either $c_2$ or $c_4$. The remaining cases yield integral circulant matrices satisfying the spectral conditions as in the following. However, none of these are realized as the concurrence matrix of a PGD.
\[\begin{array}{|c|c|c|c|c|c|}\hline
&&&\mbox{As putative } NN^T&\mbox{Feasible parameters}&\\
 \sigma & n  &r& C  & \mbox{PGD}(b, r; \alpha, \beta) & Remarks \\ \hline
 3 & 5r/3 &6l&  C[6l, l, l, l, 6l, l, l, l] &  (16l, 6l; 3l, 13l)& \mbox{Not realized} \\
 5 & r & 12l& C[12l,3l, 6l, 3l, 0, 3l, 6l, 3l] & (32l, 12l; 9l, 21l) &\mbox{Not realized} \\
 6 & 5r/6 & 48l& C[48l, 18l, 8l, 18l, 8l, 18l, 8l, 18l]& (128l, 48l; 39l, 79l) & \mbox{Undetermined} \\ \hline
\end{array}\]

\item For $k=4$, we have infinite families of circulant matrices each of which is realized as $NN^T$ for a PGD for every case as in the table below. (The case for $\sigma=1$ is missing in this table since it gives $n=4r$ and $\alpha=0$.)
\[\begin{array}{|c|c|c|c|c|c|c|}\hline
 \sigma & n  &r& C = NN^T & \mbox{PGD}(b, r; \alpha, \beta)& \mbox{Description} & \mbox{Rem.} \\ \hline
 2 & 2r & 2l & C[2l, l, 0, l, 2l, l, 0, l] & (4l, 2l; 2l, 6l) & \mbox{TD}_1(2, 2)\otimes J_{2, l} &  \ref{v=8,k=4} (1) \\
 3 & 4r/3 &3l &  C[3l, l, l, l, 3l, l, l, l] &  (6l, 3l; 4l, 8l)& \mbox{2-}(4, 2, 1) \otimes J_{2, l}& \ref{v=8,k=4} (2)  \\
 4 & r &2l & C[2l, l, l, l, 0, l, l, l] & (4l, 2l; 3l, 5l) & \mbox{TD}_l(4,2) & \ref{v=8,k=4} (3)\\
 5 & 4r/5 & 5l & C[5l, 2l, 3l, 2l, l, 2l, 3l, 2l] & (10l, 5l; 8l, 12l) & N_1\otimes J_{1, l} & \ref{v=8,k=4} (4) \\
 6 & 2r/3 & 6l & C[6l, 3l, 2l, 3l, 2l, 3l, 2l, 3l]& (12l, 6l; 10l, 14l) & N_1\otimes J_{1, l} & \ref{v=8,k=4} (5) \\ \hline
\end{array}\]
where $N_1$ is the incidence matrix of a PGD with $l=1$ and the last column is referring to Remark \ref{v=8,k=4}.

\begin{remark}\label{v=8,k=4}
We make some remarks on the five infinite families of PGDs of order $8$ whose concurrence matrices are circulant.
\begin{enumerate}
    \item[$(1)$] For $\sigma=2$, when $l=1$, the matrix $C[2,1,0,1,2,1,0,1]$ is realized as $NN^T$ for the PGD$(8, 4, 4, 2; 2, 6)$ with $\mathcal{B}=\{\{1, 2, 5, 6\}, \{2, 3, 6, 7\}, \{3, 4, 7, 8\}, \{1, 4, 5, 8\}\}$. 
    
        We can double the blocks to get a PGD$(8,8,4,4;4,12)$ which is listed in \cite{O} and obtained by the PGDS $\{0,1,4,5\} \subset \mathbb{Z}_8$.  Also, for each positive integer $l$, the PGD$(8, 4l, 4, 2l; 2l, 6l)$ is unique and can be described by a TD$_1(2,2) \otimes J_{2,l}$.  
            
    \item[$(2)$] For $\sigma = 3$, when $l=1$,  $C$ is realized as $NN^T$ for the PGD with 
    \[\mathcal{B}=\{\{1, 2, 5, 6\}, \{1, 3, 5, 7\}, \{1, 4, 5, 8\}, \{2, 3, 6, 7\}, \{2, 4, 6, 8\}, \{3, 4, 7, 8\}\}.\]
   
   With $l\ge 2$, we also obtain those PGDs listed in N9, N25, N57 in \cite{vDS}.
    
    \item[$(3)$] For $\sigma = 4$, when $l=1$, $C$ is realized as $NN^T$ for a PGD$(8, 4, 4, 2;3, 5)$ which is isomorphic to a TD$_1(4, 2)$. Furthermore,  $C$ is realized as $NN^T$ for a TD$_2(4, 2)$ when $l=2$ and a TD$_3(4,2)$ when $l=3$, as well. The block sets of these two PGDs found by computer search are, respectively,\\
   \indent \qquad \quad $ \left \{\begin{array}{cccc} \{1, 2, 3, 4\}, & \{1, 2, 4, 7\}, & \{1, 3, 6, 8\}, & \{1, 6, 7, 8\},\\
     \{2, 3, 5, 8\}, & \{2, 5, 7, 8\}, &  \{3, 4, 5, 6\}, &\{4, 5, 6, 7\}\end{array}\right \}, \ \mbox{ and }$  \\
     $\left \{\begin{array}{cccccc} \{1, 2, 3, 4\}, & \{1, 2, 4, 7\}, & \{1, 2, 7, 8\}, & \{1, 3, 4, 6\},& \{1, 3, 6, 8\}, &
     \{1, 6, 7, 8\},\\  \{2,3, 4, 5\}, &  \{2, 3, 5, 8\}, &\{2, 5, 7, 8\}, & 
     \{3, 5, 6, 8\}, & \{4, 5, 6, 8\}, & \{4, 5, 6, 8\}\end{array}\right \}$.
    
    From these two, for any $l \geq 4$, the block set of a TD$_l(4,2)$ is obtained, and $C[2l, l, l, l, 0, l, l, l]$ is realized as $NN^T$ for the TD$_l(4,2)$  (cf. Ex. \ref{PGDS1}). Those with $l=2, 3, 4, 5, 6$ appear as N2, N8, N16, N35, N56, resp., in \cite{vDS}. 
    
    \item[$(4)$] For $\sigma = 5$, $C[5, 2, 3, 2, 1, 2, 3, 2]$ is realized as $NN^T$ for the PGD$(8, 10, 4, 5; 8,12)$ with $\mathcal{B}=\left \{\begin{array}{ccccc} \{1, 2, 3, 4\}, & \{1, 2, 3, 8\}, & \{1, 3, 5, 7\}, & \{1, 4, 6, 7\},&
     \{1, 6, 7, 8\},\\  \{2, 4, 5, 7\}, &  \{2, 4, 6, 8\}, &\{2, 5, 7, 8\}, & \{3, 4, 5, 6\}, & \{3, 5, 6, 8\}\end{array}\right \}$.\\ (The incidence matrix for this design is found via a computer search.) 
    This PGD and the PGD$(8, 20, 4, 10; 16, 24)$ are listed as N4 and N34, resp., in \cite{vDS}. 
    
    \item[$(5)$] For $\sigma=6$, our computer search confirms that $C[6, 3, 2, 3, 2,3, 2, 3]$ is realized as $NN^T$ for a PGD$(8, 12, 4, 6;10, 14)$. (This is N7 in \cite{vDS}). The block set $\mathcal{B}$ for this PGD found via computer search is given by\\
    $\left \{\begin{array}{cccccc} \{1, 2, 3, 4\}, & \{1, 2, 3, 8\}, & \{1, 2, 6, 7\}, & \{1, 4, 7, 8\},& \{1, 4, 5, 6\}, & \{1, 5, 6, 8\}\\
     \{2, 3, 5, 6\}, & \{2, 4, 5, 7\},  &\{2, 5, 7, 8\}, & \{3, 4, 6, 7\}, & \{3, 4, 5, 8\},& \{3, 6,7, 8\}\end{array}\right \}$.

 For $l=2$, $C[12, 6, 4, 6, 4, 6, 4]$ is $NN^T$ for a PGD$(8, 24, 4, 12; 20, 28)$, one of 56 listed as N55 in \cite{vDS}.
\end{enumerate}
\end{remark}
\end{enumerate}

We have found a PGD for every parameter set that is also in \cite{vDS}.  In some cases, we have found all the PGDs in \cite{vDS}.  In other cases, we are missing some. The cases for $k=5, 6$ are omitted as all PGDs with $k=5$ are the complementary designs for those with $k=3$. The PGDs with $k=2$ (complementary of $k=6$) are improper and omitted. 
This concludes our analysis and classification of PGDs of order $8$ with circulant concurrence matrices. 


\subsection{PGDs of order 9}{\label{A9}}

The eigenvalues of the $9\times 9$ circulant matrix $C=C[c_0, c_1, c_2, c_3, c_4, c_4, c_3, c_3, c_1]$ are given by, for $j=0, 1, 2,\dots, 8$, 
$$\theta_j = c_0 + 2c_1\cos\left(\frac{2j\pi}{9}\right) + 2c_2\cos\left(\frac{4j\pi}{9}\right) + 2c_3\cos\left(\frac{2j\pi}{3}\right) + 2c_4\cos\left(\frac{8j\pi}{9}\right).$$ 
Namely, 
\[\begin{array}{ll}
f(\omega^1) =f(\omega^8) = & c_0 + 2c_1\cos(2\pi/9) +2c_2\cos(4\pi/9) -c_3 + 2c_4 \cos(8\pi/9),\\
f(\omega^2) =f(\omega^7) = & c_0 + 2c_1\cos(4\pi/9) +2c_2\cos(8\pi/9) -c_3 + 2c_4 \cos(2\pi/9),\\
f(\omega^4) =f(\omega^5) = & c_0 + 2c_1\cos(8\pi/9) +2c_2\cos(2\pi/9) -c_3 + 2c_4 \cos(4\pi/9),\\
f(\omega^3) =f(\omega^6) = &  c_0 \ \ - \ \ c_1 \ \ -\ \  c_2 \ \ +\ \  2c_3 \ \  -\ \ c_4 .\end{array}\]

For $C$ to be realized as the concurrence matrix of a PGD$(9, b, k, r;\alpha, \beta, n)$, $C$ must have integral eigenvalues $kr, n$ and $0$ with multiplicities $1, \sigma$ and $8-\sigma$, respectively. (Note that   $n\sigma = r(9-k)$, together with $b=9r/k$, 
$\alpha=k(kr-n)/9$ and $\beta=n+k(kr-n)/v$.)
Under these constraints and the fact that $\cos\left(\frac{2\pi}{9}\right) + \cos\left(\frac{4\pi}{9}\right) + \cos\left(\frac{8\pi}{9}\right) = 0$, the feasible integral values of the $c_i$'s must satisfy $c_1=c_2=c_4$. Then 
$\theta_1 = \theta_2 = \theta_4 = \theta_5 = \theta_7 = \theta_8 = c_0-c_3$, $\theta_3=\theta_6=c_0 -c_1-c_2+2c_3-c_4$, and we have

\begin{tabular}{ |c||c|c|c|c| } 
 \hline
 $j$ & $\theta_j$ &  $\theta_j$ expected & multiplicity & remark\\
 \hline
 $0$ & $c_0 + 2c_1 + 2c_2 + 2c_3 + 2c_4$ & $kr$ & 1 & \\ 
 $1,2,4,5,7,8$ & $c_0 - c_3$ & either $0$ or $n$ & 6 & $c_1 = c_2=c_4$\\ 
 $3$, $6$ & $c_0-c_1-c_2+2c_3-c_4$ & either $0$ or $n$ & 2  & \\ 
 \hline
\end{tabular}

Suppose that Spec$(C)=[kr^1, n^{\sigma}, 0^{8-\sigma}] \ \mbox{with } \ \sigma = \frac{r(9-k)}{n}$. Then the only possible value for $\sigma$ is either 2 or 6. 
\[\begin{array}{|c|c||c|c|c||c|} \hline
 \sigma & n &\theta_0(=kr) &  \theta_1, \theta_2, \theta_4, \theta_5, \theta_7& \theta_3, \theta_6 & \mbox{feasible }(c_0,\  c_1,\  c_2, \ c_3,\  c_4) \\ \hline
  2 &\frac12r(9-k) &9r-2n&  0 & n & (r, \ r-\frac{n}{3}, \ r-\frac{n}{3}, \ r, \ r-\frac{n}{3})   \\ \hline
 6 & \frac16r(9-k) &9r-6n&  n & 0 & (r, r-\frac{2n}{3}, r-\frac{2n}{3},  r-n, r-\frac{2n}{3})  \\ \hline
\end{array}\]

For $k=3$ and $k=4$, the feasible solutions give the following putative concurrence matrices and corresponding parameters of PGDs.
\[\begin{array}{|c||crlr||crlr|} \hline
 & &  k \quad = & \quad 3 &  & &k \quad =  & 4  &\\ 
\sigma & n & [c_0, c_1, c_2, c_3, c_4] & (b, r; \alpha,\beta) & \mbox{Rem.} &
n &   [r, c_1, c_2, c_3, c_4] & (b, \ r; \ \alpha, \ \beta) & \mbox{Rem.}\\ \hline 
2 & 3r &  [r, \  0,\   0,\  r, \  0] & (3r, r; 0, 3r) & \ref{v=9k=3} (1)& 
\frac{5r}{2} & [r, \ \frac{r}{6}, \ \frac{r}{6},\  r, \ \frac{r}{6}]  &  (\frac{9r}{4}, r; \frac{2r}{3}, \frac{19r}{6}) & \ref{v=9k=3}(3)\\
6 &  3l& [3l,\ l,\ l, \ 0, \ l]&  (9l,  3l; 2l, 5l)&  \ref{v=9k=3}(2) &
  \frac{5r}{6}& [r, \frac{4r}{9}, \frac{4r}{9}, \ \frac{r}{6},  \frac{4r}{9}]  & (\frac{9r}{4}, r;  \frac{38r}{27}, \frac{121r}{54}) &  \ref{v=9k=3}(3) \\ \hline
\end{array}\]

\begin{remark}\label{v=9k=3}
\begin{enumerate}
\item[$(1)$] For $k=3$ and $\sigma = 2$, the putative circulant matrix $C[r, 0, 0, r, 0, 0, r, 0, 0]$ is realized as the concurrence matrix of an improper PGD$(9,3r,3,r;0,3r)$ with multi-block set $\mathcal{B} = \{ \{1,4,7\}^r, \{2,5,8\}^r, \{3,6,9\}^r\}$.

\item[$(2)$] For $k=3$ and $\sigma = 6$, $C[3l, l, l, 0, l, l, 0, l, l]$ is realized  as the concurrence matrix of a PGD$(9, 9l, 3, 3l; 2l, 5l)$ for all $l\in\mathbb{N}$. We describe the first three cases of $l$:
\begin{enumerate}
    \item For $l=1$ and $n=r=3$, $C[3, 1, 1, 0, 1, 1, 0, 1, 1]$ is  $NN^T$ for the PGD$(9,9,3,3;2,5)$ which is isomorphic to N5 in \cite{vDS}, a TD$_1(3,3)$ and a PG$(3, 2, 2)$ in Example \ref{PG}. Here is the block set of this PGD:
\[\{\{1, 2, 3\}, \{1, 5, 6\}, \{1, 8, 9\}, \{2, 4, 9\}, \{2, 6, 7\}, \{3, 4, 5\},\{3, 7, 8\}, \{4, 6, 8\}, \{5, 7, 9\}\}\]

    \item For $l=2$ and $n=r=6$, $C$ is realized as $NN^T$ for a PGD$(9,18,3,6;4,10)$ which is N29 in \cite{vDS}.  We can take the incidence matrix of a TD$_2(3,3)$. However, as it is indicated in \cite{vDS}, there are four PGDs with the same parameters. Here is one of the block sets we have from our computer search:
\[\left \{\begin{array}{l} \{1, 2, 3\}, \{1, 2, 6\}, \{1, 3, 8\}, \{1, 5, 6\}, \{1, 5, 9\}, \{1, 8, 9\},\{2, 3, 7\}, \{2,4, 6\}, \{2, 4, 9\}\\
\{2, 7, 9\}, \{3, 4,5\}, \{3, 4, 8\}, \{3, 5, 7\}, \{4, 5, 9\}, \{4, 6, 8\},\{5, 6, 7\}, \{6, 7, 8\}, \{7, 8, 9\}\end{array}\right \}\]

   \item For $l=3$ and $n=r=9$, the corresponding circulant matrix is realized as the concurrence matrix of a PGD$(9,27,3,9;6,15)$ in three distinct ways. 
   \begin{enumerate}
   \item The first way is $3.$TD$_1(3,3)$ by taking $N \otimes J_{1,3}$ where $N$ is the matrix of TD$_1(3,3)$ above.  
   \item The second way is by taking the union (as a multi-set) of blocks from the PGDs when $n=r=3$ and $n=r=6$ (i.e. concatenating the two incidence matrices in parts (a) and (b) which are described as  TD$_1(3,3) \ \uplus $ TD$_2(3,3)$, where `$\uplus$' denotes the multi-set union; so $9$ blocks of which are not repeated and the other $9$ which are repeated exactly once.  
   \item The third way is by taking the incidence matrix of a TD$_3(3,3)$ that does not have repeated blocks.   
   \end{enumerate}
  
\item In general, for every $l\ge 1$, the circulant matrix $C[3l, l, l, 0, l, l, 0, 1, 1]$ becomes the concurrence matrix of a PGD$(9, 9l, 3, 3l; 2l, 5l)$ which is isomorphic to $\biguplus\limits_{i=1}^3 n_i.\mbox{TD}_i(3, 3)$ over all possible nonnegative integer triples $(n_1, n_2, n_3)$ with $l=n_1+2n_2+3n_3$.

\end{enumerate}

\item[$(3)$] For $k=4$, whether $\sigma = 2$ or 6, it is shown that none of the putative circulant matrices are realized as the concurrence matrix of a PGD.
\end{enumerate}
\end{remark}
\subsection{PGDs of order 10}{\label{A10}}


Similar analysis on a $10\times 10$ symmetric circulant  matrix $C=C[c_0, c_1, c_2, c_3, c_4, c_5]$ gives us two infinite families of symmetric circulant matrices each of which is realized as the concurrence matrix for a PGD for the following combinations of $k$ and $\sigma$.
\[\begin{array}{|c|c|c|c|c|c|}\hline
 k & \sigma & n & \mbox{Feasible } C & \mbox{PGD}(10, b, k, r;\alpha, \beta) & \mbox{Rem.}\\ \hline
 4 & 4 & 6l (l\ge 1)& C[4l,l,l,l,l, 4l]  & (10, 10l, 4, 4l; 4l, 10l) & (1)\\
5 & 5 & 2l (l\ge 3) & C[2l, l, l, l, l, 0] & (10, 4l, 5, 2l; 4l, 6l) & (2)\\ \hline
\end{array}\]

\begin{remark}
\begin{enumerate}
\item[$(1)$] For $k=4$ and $\sigma=4$, the circulant matrices are the concurrence matrices of the following PGDs: 
\begin{enumerate}
    \item For $l=1$, PGD($10,10,4,4;4,10)$ which is N10 in \cite{vDS};
    \item for $l=2$,  PGD($10,20,4,8;8,20)$ which is N45 in \cite{vDS},
\end{enumerate}
 Furthermore, for each $l\ge 1$, this PGD is unique (can easily be shown by hand).
\item[$(2)$] When $k=5$ and $\sigma=5$, the circulant matrix is realized as the concurrence matrix of a PGD isomorphic to TD$_l(5, 2)$. 
  By incrementing $l$, we get the following PGDs:
\begin{enumerate}
    \item PGD($10,8,5,4;8,12)$ which is not in \cite{vDS} and the corresponding $C$ is realized by the incidence matrix $N_1$ on the left below (which is also that of a TD$_2(5,2)$). 
    \item PGD($10,12,5,6;12,18)$ which is N14 in \cite{vDS} and $C=NN^T$ can be obtained by $N_2$ on the right below.
    \item PGD($10,16,5,8;16,24)$ which is N26 in \cite{vDS},
    \item PGD($10,20,5,10;20,30)$ which is N46 in \cite{vDS},
    \item PGD($10,24,5,12;24,36)$ which is N77 in \cite{vDS}.
    {\small 
\[N_1 = \left [\begin{array}{cccccccc}
1 & 1 & 1 & 1 & 0 & 0 & 0 & 0\\
1 & 1 & 0 & 0 & 1 & 1 & 0 & 0\\
1 & 1 & 0 & 0 & 0 & 0 & 1 & 1\\
1 & 0 & 1 & 0 & 1 & 0 & 1 & 0\\
1 & 0 & 1 & 0 & 0 & 1 & 0 & 1\\
0 & 0 & 0 & 0 & 1 & 1 & 1 & 1\\
0 & 0 & 1 & 1 & 0 & 0 & 1 & 1\\
0 & 0 & 1 & 1 & 1 & 1 & 0 & 0\\
0 & 1 & 0 & 1 & 0 & 1 & 0 & 1\\
0 & 1 & 0 & 1 & 1 & 0 & 1 & 0
\end{array}\right ] \quad
N_2 = \left [\begin{array}{cccccccccccc} 
1 & 1 & 1 & 1 & 1 & 1 & 0 & 0 & 0 & 0 & 0 & 0\\
1 & 1 & 1 & 0 & 0 & 0 & 1 & 1 & 1 & 0 & 0 & 0\\
1 & 1 & 0 & 1 & 0 & 0 & 1 & 0 & 0 & 1 & 1 & 0\\
1 & 1 & 0 & 0 & 1 & 0 & 0 & 1 & 0 & 1 & 0 & 1\\
1 & 1 & 0 & 0 & 0 & 1 & 0 & 0 & 1 & 0 & 1 & 1\\
0 & 0 & 0 & 0 & 0 & 0 & 1 & 1 & 1 & 1 & 1 & 1\\
0 & 0 & 0 & 1 & 1 & 1 & 0 & 0 & 0 & 1 & 1 & 1\\
0 & 0 & 1 & 0 & 1 & 1 & 0 & 1 & 1 & 0 & 0 & 1\\
0 & 0 & 1 & 1 & 0 & 1 & 1 & 0 & 1 & 0 & 1 & 0\\
0 & 0 & 1 & 1 & 1 & 0 & 1 & 1 & 0 & 1 & 0 & 0
\end{array}\right ]\]}
\end{enumerate}
For each $l\ge 4$, the associated PGD is isomorphic to TD$_l(5,2)$, and can be attained from the above TD$_2(5,2)$ and TD$_3(5,2)$ as described before.

\item[$(4)$]  We note that among the PGDs of order 10, the ones labeled as N23 (the PGD with parameters $(10, 15, 4, 6; 10, 15)$), N36 (the PGDs with parameters $(10, 18, 5, 9; 18, 27)$) and N86 (the PGD$(10, 25, 4, 10; 10, 25)$)  in \cite{vDS} are not attained from our analysis of concurrence matrices. 

\end{enumerate}
\end{remark}


\subsection{PGDs of order 12}
For the order 12 case, we find out that there are infinite families of symmetric circulant matrices that are realized as the concurrence matrices of PGDs with block sizes 3, 4 and 6. We summarize our results here with some key descriptions of the basic examples for important cases, and for the further details we refer the readers to \cite{Tr}.

A symmetric circulant matrix $C = C[c_0, c_1, c_2, c_3, c_4, c_5, c_6, c_5, c_4, c_3, c_2, c_1]$ is feasible to be realized as the concurrence matrix, $NN^T$, of a PGD$\left(12, \frac{12r}{k}, k, r; \frac{k(kr-n)}{12}, n+\frac{(kr-n)}{12}\right)$ when
$$\text{Spec}(C)=\text{Spec}(NN^T) = [kr^1, n^{\sigma}, 0^{11-\sigma}] \,\,\, \text{where} \,\,\, \sigma = \frac{r(12-k)}{n}$$
with the eigenvalues of $C$ are given by, for $j=0, 1, \dots, 11$,  {\small $$\theta_j  = c_0 + 2c_1\cos\left(\frac{j\pi}{6}\right) + 2c_2\cos\left(\frac{j\pi}{3}\right) + 2c_3\cos\left(\frac{j\pi}{2}\right) + 2c_4\cos\left(\frac{2j\pi}{3}\right) +2 c_5\cos \left (\frac{5j\pi}{6}\right )+c_6 \cos\left(j \pi\right).$$}

Using the above conditions on the spectra of $NN^T$ along with the integrality conditions on the eigenvalues of a circulant matrix, we identify all feasible symmetric circulant matrices and corresponding PGDs of order 12. 

\begin{theorem} There are at least 11 infinite families of $12\times 12$ symmetric integral circulant matrices that are realized as concurrent matrices of PGDs. All these matrices and the corresponding PGDs are listed in the following table.
\[\begin{array}{|l|c|c|rrrr|c|cc|} \hline
\mbox{Symmetric circulant}& &&  \mbox{PGD}& &&&&&\\
C[c_0, c_1, c_2, c_3, c_4, c_5, c_6] & k & n^{\sigma} & b, &  r ; & \alpha, & \beta   &  \mbox{Description} & \# & \mbox{Remark}  \\ \hline
C[4l,\ l,\ l,\ 0,\ l,\ l, \ 0] & 3 & (4l)^9 & 16l,  & 4l; & 2l, & 6l &   \biguplus_{i=1}^{4} n_i.\mbox{TD}_i(3,4) &\ge1  &R(1)\\
\hline
C[3l, \ l, \ 0,\ l,\ 0,\ l, \ 3l] & 4 & (6l)^4 & 9l, & 3l; & 2l, &8l &  \mbox{TD}_1(2,3) \otimes J_{2,l} & \exists ! &R(2(a)) \\

C[5l, \ l,\ l,\ l,\ l,\ l,\ 5l] & 4 & (8l)^5 & 15l,  & 5l; & 4l, &12l &   2-(6,2,1) \otimes J_{2,l} & \exists ! &R(2(b)) \\

C[4l, \ l,\ l,\ 2l,\ l, \ l,\ 0] & 4 & (4l)^8 & 12l,  & 4l; & 4l, &8l &  computer & \ge 1 & R(2(c)i)\\

C[3l,\ l,\ l,\ l,\ 0,\ l,\ l] & 4 & (3l)^8 & 9l,  & 3l; & 3l, &6l &  \mbox{TD}_1(4,3) \otimes J_{1,l} & \ge 1 & R(2(c)ii) \\
\hline
C[2l,\ l,\ l,\ 0,\ l,\ l, \ 2l] & 6 & (4l)^3 & 4l,  & 2l; & 4l, &8l &  D_{1,6,2l} \otimes J_{2,1} &\ge 1  &R(3(a)i) \\

C[3l, \ l,\ l,\ l,\ 3l,\ l,\ l] & 6 & (6l)^3 & 6l,  & 3l; & 6l, &12l &  2-(4,2,1) \otimes J_{3,l} & \exists ! & R(3(a)ii)\\

C[5l, 2l, 2l, 2l, 2l, 2l,  5l] & 6 & (6l)^5 & 10l , & 5l; & 12l, &18l &  2-(6,3,2) \otimes J_{2,l} & \exists ! &R(3(b)) \\

C[2l,\ l,\ l,\ l,\ l,\ l,\ 0] & 6 & (2l)^6 & 4l , & 2l; & 5l, &7l &  \mbox{TD}_l(6,2), l\ge 2 & \ge 1 &R(3(c)) \\

C[7l, 3l, 4l, 3l, 4l, 3l,  l] & 6 & (6l)^7 & 14l,  & 7l; & 18l, &24l &  computer & \ge 1 & R(3(d))\\

C[6l, 3l, 3l, 2l, 3l, 3l, 2l] & 6 & (4l)^9 & 12l,  & 6l; & 16l, &20l &  computer & \ge 1 &R(3(e)) \\
\hline
\end{array}\]
\end{theorem}
\begin{proof} We give some narrative comments on the results in the following remark instead of providing a tedious derivation of the entries in the table. For the details, we refer to \cite{Tr}.
\end{proof}

\begin{remark}
\begin{enumerate}
\item[(1)] For $k=3$, the only feasible case occurs when $\sigma=9$.  
\begin{enumerate}
\item With $l=1$, $C[4, 1, 1, 0, 1, 1, 0]$  is realized as $NN^T$ for a PGD($12,16,3,4;2,6)$ which is isomorphic to a TD$_1(3,4)$.  This PGD is listed as N37 in \cite{vDS}. 
\item With $l=2$, $C[8, 2, 2, 0, 2, 2]$ is realized as $NN^T$ for two PGDs, TD$_2(3,4)$ and  2.TD$_1(3,4)$ ($=$TD$_1(3,4) \otimes J_{1,2}$), with the same parameters $(12,32,3,8;4,12)$.  
\item With $l=3$, $C[12, 3, 3, 0, 3, 3]$ is realized as $NN^T$ for three PGDs with the same parameters $(12,48,3,12;6,18)$; namely,  \begin{enumerate}
    \item TD$_3(3,4)$,
    \item 3.TD$_1(3, 4)$, with its incidence matrix $N \otimes J_{1,3}$ (where $N$ is the incidence matrix of  TD$_1(3,4)$),
    \item TD$_1(3,4)\uplus$ TD$_2(3,4)$.
\end{enumerate}
\item With $l\ge 4$, similarly, $C[4l, l, l, 0, l, l,0]$ are realized as the concurrence matrices for PGDs
with parameters $(12, 16l, 3, 4l; 2l, 6l)$ belonging to the family
\[\left \{ \biguplus_{i=1}^4 n_i.TD_i(3,4) \ :\ \mbox{for all nonnegative solutions } (n_1, n_2, n_3, n_4) \mbox{ for }  l= \sum\limits_{i=1}^4 in_i\right \}.\]
\end{enumerate}

\item[(2)] 
For $k=4$, there are three different values 4, 5, 8 of $\sigma$ which give rise to four different infinite families of PGDs.  
\begin{enumerate}
    \item ($\sigma=4$): For $l=1$, $C[3, 1, 0, 1, 0, 1, 3]$ is $NN^T$ for a unique PGD $TD_1(2,3) \otimes J_{2,1}$ with parameters $(12, 9, 4, 3;2, 8)$. Here TD$_1(2,3)$ is the PGD coming from bipartite graph $K_{3,3}$ with point set $\{1, 3, 5\}\cup\{2, 4, 6\}$ and edge set as the block set. This PGD is not listed in \cite{vDS} as it is considered as  improper there. For $l=2$, $C[6, 2, 0, 2,0, 2, 6]$ is realized as $NN^T$ for a PGD$(12, 18, 4, 6; 4, 16)$ which is described as a $TD_1(2,3)\otimes J_{2,2}$ and is listed as N48 in \cite{vDS}.
    
        \item ($\sigma=5$): For each $l$, $C[5l, l, l, l, l, l, 5l]$ gives rise to a unique PGD.  For $l=1$, the PGD is isomorphic to $D\otimes J_{2, 1}$ where $D$ is the 2-$(6, 2, 1)$ design, which is listed as   N33 in \cite{vDS}.
        
        \item ($\sigma=8$): 
    \begin{enumerate}
        \item $C[4, 1, 1, 2, 1, 1, 0]$ is $NN^T$ for the PGD$(12, 12, 4, 4; 4, 8)$ which is N19 in \cite{vDS}. We have the following block set  with $P=\{0, 1, 2, \dots, 9, a, b\}$ via computer search:\\
$\mathcal{B}=\left \{ \begin{array}{c} \{0, 1, 3, 4\}, \{0, 2, 3, 5\}, \{0, 7, 9, a\},   \{0, 8, 9, b\}, \{1, 2, a, b\}, \{1, 4, 6, 9\}, \\ 
\{1, 5, 8, a\}, \{2, 4, 7, b\}, \{2, 5, 6, 9\}, \{3, 6, 7, a\}, \{3, 6, 8, b\}, \{4, 5, 7, 8\} \end{array}  \right \}$.    
        
          By multipling the block set, we get a PGD for each member of this infinite family.

        \item With $l=1$, $C[3, 1, 1, 1, 0, 1, 1]$ is $NN^T$ for a TD$_1(4,3)$, a PGD$(12,9, 4, 3; 3, 6)$.\\ For each $l\ge 2$, $C[3l, l, l, l, 0, l, l]$ is realized as $NN^T$ for the PGD $l.$TD$_1(4,3)$.
        (The case for $l=2$ is N47 in \cite{vDS}.)
    \end{enumerate}
\end{enumerate}

\item[(3)] 
For $k=6$, there are five different values, 3, 5, 6, 7, 9 of $\sigma$ which give rise to $6$ different infinite families of PGDs.  
\begin{enumerate}
    \item ($\sigma=3$):
    \begin{enumerate}
        \item For $l=1$, $C[2, 1, 1, 0, 1, 1, 2]$ is $NN^T$ for the PGD $\mathcal{D}=TD_1(3,2)\otimes J_{2,1}$. For each $l\ge 2$, $C[2l, l, l, 0, l, l, 2l]$ is $ NN^T$ for $l.\mathcal{D}=\mathcal{D}\otimes J_{1, l}$.
        The PGDs $3.\mathcal{D}$, $4.\mathcal{D}$ and $5.\mathcal{D}$ are listed as N22, N41 and N62, respectively,  in \cite{vDS}.
        \item For $l=1$, $C[3, 1, 1, 1, 3, 1, 1]$ is $NN^T$ for $\mathcal{D}= D\otimes J_{3,1}$ where $D$ is a 2-$(4, 2,1)$.
        For each $l\ge 2$, $C[3l, l, l, l, 3l, l, l]$ is $NN^T$ for $l.\mathcal{D}=D\otimes J_{3,l}$ with parameters $(12, 6l, 6, 3l; 6l, 12l)$.
        Among those, the cases for $l=2$ and 3 are listed as N22 and N50, respectively, in \cite{vDS}.
    \end{enumerate}
    \underline{Note}: We have found two PGDs with the same parameters but with different circulant concurrence matrices, namely, the two examples that are N22 in \cite{vDS}.
    
    \item ($\sigma=5$): For $l=1$, $C[5, 2, 2, 2, 2, 2, 5]$ is realized as $NN^T$ for $D\otimes J_{2,1}$, the PGD$(12, 10, 6, 5; 12, 18)$,  where $D$ is a 2-$(6, 3,2)$ design with block set \[\left \{\begin{array}{c} \{1, 2, 3\}, \{1, 2, 4\}, \{1, 3, 5\}, \{1, 4, 6\}, \{1, 5, 6\},\\ \{2, 3, 6\}, \{2, 4, 5\}, \{2, 5, 6\}, \{3, 4, 5\}, \{3, 4, 6\}\end{array}\right \}.\]
    With $l\ge 2$, $C[5l, 2l, 2l, 2l, 2l, 2l, 5l]$ is $NN^T$ for $D\otimes J_{2,l}$, and  the case when $l=2$ is identified with N61 in \cite{vDS}.
    
    \item ($\sigma=6$): For $l=1$, $C[2, 1, 1, 1, 1,1, 0]$ is not realized. For $l=2$, $C[4, 2, 2, 2, 2, 2, 0]$ is realized as $NN^T$ for a TD$_2(6,2)$, a PGD$(12, 8, 6, 4; 10, 14)$ whose block set is 
    \[\mathcal{B}_2(6,2)=
    \left \{\begin{array}{c} \{0, 1,2, 3, 4, 5\},  \{0, 1, 2, 9, a, b\},  \{0, 3, 4, 7, 8, b\},  \{0, 5, 7, 8, 9, a\}\\
     \{1, 3, 5, 6, 8, a\},  \{1, 4, 6, 8, 9, b\},  \{2, 3, 6, 7, a, b\},  \{2, 4, 5, 6, 7, 9\}\\ \end{array}\right \}.\]
     For $l=3$, $C[6, 3, 3, 3, 3, 3,0]$ is $NN^T$ for a TD$_3(6,2)$, a PGD$(12, 12, 6, 6, 15, 21)$  whose block set is
    \[\mathcal{B}_3(6,2)=
    \left \{\begin{array}{l} \{0, 1,2, 3, 4, 5\},  \{0, 1, 2,3, 4, b\},  \{0, 1, 5, 8, 9,a\},  \{0,2, 5, 7, 9, a\},\\
    \{0, 3, 7, 8, a, b\}, \{0, 4, 7, 8, 9, b\}, \{1, 2, 6, 9, a, b\}, \{1, 3, 5, 6, 8, a\}\\
     \{1, 4, 6, 8,9, b\},   \{2, 3, 6, 7, a, b\},  \{2, 4, 5, 6, 7, 9\}, \{3, 4, 5, 6, 7, 8\}\\ \end{array}\right \}.\]
      
      For each $l\ge 4$, $C[2l, l, l, l, l, l, 0]$ is realized as $NN^T$ for 
      a TD$_l(6,2)$ which can be attained as a suitable combination of $\mathcal{B}_2(6,2)$ and $\mathcal{B}_3(6,2)$.
      The cases $l=3$ and 4 are listed as N21 and N39, respectively in \cite{vDS}.
    
    \item ($\sigma=7$): For $l=1$, $C[7, 3, 4, 3, 4, 3, 1]$ is realized as $NN^T$ for  a PGD$(12, 14, 6, 7; 18, 24)$ which is listed as N27 in \cite{vDS}, but we do not know if it is unique. Via computer search we have the following block set for this case:
     \[\left \{\begin{array}{l} \{0, 1,2, 3, 4, 5\},  \{0, 1, 2,3, 4, b\},  \{0, 1, 5, 8, 9,a\},  \{0,2, 4, 6, 8, a\}, \{0, 2, 5, 7, 9, a\}\\
    \{0, 3, 7, 8, a, b\}, \{0, 4, 7, 8, 9, b\}, \{1, 2, 6, 9, a, b\}, \{1, 3, 5, 6, 8, a\}, \{1, 3, 5, 7, 9, b\}, \\
     \{1, 4, 6, 8,9, b\},   \{2, 3, 6, 7, a, b\},  \{2, 4, 5, 6, 7, 9\}, \{3, 4, 5, 6, 7, 8\}\\ \end{array}\right \}.\]
    
    \item ($\sigma=9$): For $l=1$, $C[6, 3, 3, 2, 3, 3, 2]$ is realized as $NN^T$ for a PGD$(12, 12, 6, 6; 16, 20)$ which is listed as N20 in \cite{vDS}. We find the following block set for this design via computer search:
    \[\left \{\begin{array}{l} \{0, 1,2, 3, 4, 8\},  \{0, 1, 2,6, 7, b\},  \{0, 1, 2, 5, 9,a\},  \{0,3, 7, 8, a, b\},\\
    \{0, 4, 5, 7, 9, b\}, \{0, 4,5, 6,  8, a\}, \{1, 3,4,5, 6, b\}, \{1, 3, 5, 7, 8, 9\}\\
     \{1, 6, 8,9, a, b\},   \{2, 3, 4, 9, a, b\},  \{2, 3, 5, 6, 7, a\}, \{2, 4, 6, 7, 8,9\}\\ \end{array}\right \}.\]
      
\end{enumerate}

\item[(4)] For $k=6$ and  $\sigma=8$, there is a putative circulant matrix $C[8l, 4l, 4l, 4l, l, 4l, 4l]$ which may be realized as $NN^T$ for a PGD$(12, 16l, 6, 8l; 21l, 27l)$. Although N38 in \cite{vDS} indicates that for $l=1$, it is realized as the concurrence of PGD$(12, 16, 6, 8; 21, 27)$ the size is too large for us to determine this case.
\end{enumerate}

\end{remark}


\section{Concluding remarks}

We notice that there are some PGDs that are listed in \cite{vDS} that our method did not find.  

The following PGDs listed in \cite{vDS} were not found using our method:
\begin{enumerate}
    \item Order 8: We have found a PGD for every parameter set that is listed in \cite{vDS}.  For some parameter sets, we have found all the PGDs that were listed, in other cases, we are missing some, and yet in other cases, we do not know if we can find all the listed PGDs. 
    \item Order 9: We have found the PGD that is N4 and one of the $4$ PGDs that are listed in N29 in \cite{vDS}.  But, our computer search did not realize any of the following listed in \cite{vDS}:
    \begin{enumerate}
        \item N12, a $\overline{\mbox{TD}_1(2,3)} \otimes J_{1,2}$ 
        \item N30, a $\overline{\mbox{TD}_1(2,3)} \otimes J_{1,3}$
        \item N71, a $\overline{\mbox{TD}_1(2,3)} \otimes J_{1,4}$
    \end{enumerate}
    \item Order 10: We have found a PGD for every parameter set that is also listed in \cite{vDS} except for the following three.  
    \begin{enumerate}
        \item N23, a $\overline{\mbox{D}} \otimes J_{1,3}$ where $D$ is a 2-$(5, 2,1)$ design
        \item N36, a $\overline{\mbox{D}} \otimes J_{1,3}$ where $D$ is a $2$-$(6,3,2)$ and via computer search
        \item N86, a $\overline{\mbox{D}} \otimes J_{1,5}$ where $D$ is a 2-$(5, 2,1)$ design
    \end{enumerate}
    \item Order 12:  We have found a PGD for every parameter set that is also listed in \cite{vDS} except for the following three.  
    \begin{enumerate}
        \item N40, a $\overline{\mbox{TD}_3(4,2)} \otimes J_{1,2}$, and via computer search
        \item N38 \& N49.  The parameters for the PGDs listed here are too large for a computer search to find a PGD.  Thus, for these two, we cannot determine if a PGD can be found using our method.
    \end{enumerate}
\end{enumerate}


\end{document}